\def\SUBMITTED{}
\documentclass[journal,onecolumn,11pt]{IEEEtran}
\usepackage{amsmath,amssymb,amsthm}
\usepackage{array,booktabs}
\usepackage{url}
\usepackage[hidelinks]{hyperref}
\hypersetup{pdftitle={Optimal Reduced Representatives of Quadratic APN Permutations in Odd Dimension},pdfauthor={Hiroki Fukui}}
\newtheorem{theorem}{Theorem}[section]
\newtheorem{lemma}[theorem]{Lemma}
\newtheorem{proposition}[theorem]{Proposition}
\newtheorem{corollary}[theorem]{Corollary}
\newtheorem*{mainthm}{Main Theorem}
\theoremstyle{definition}
\newtheorem{definition}[theorem]{Definition}
\newtheorem{remark}[theorem]{Remark}
\newtheorem{example}[theorem]{Example}
\newcommand{\F}{\mathbb F}
\newcommand{\Tr}{\mathrm{Tr}}
\newcommand{\red}{\mathrm{red}}
\newcommand{\GR}{\mathrm{GR}}
\newcommand{\GL}{\mathrm{GL}}
\newcommand{\rk}{\operatorname{rk}}
\newcommand{\rank}{\operatorname{rank}}
\newcommand{\row}{\operatorname{row}}
\newcommand{\im}{\operatorname{im}}
\newcommand{\Pf}{\operatorname{Pf}}
\newcommand{\bG}{\mathbf G}
\newcommand{\lean}[1]{\texttt{#1}}

\newcommand{\venuetitle}{Optimal Reduced Representatives of Quadratic APN Permutations in Odd Dimension}
\begin{document}
\sloppy\emergencystretch=3em
\title{\venuetitle}
\author{Hiroki~Fukui%
\thanks{H. Fukui is with the Research Institute of Criminal Psychiatry / Sex Offender Medical Center, Tokyo, Japan, and with the Department of Neuropsychiatry, Kyoto University, Kyoto, Japan (e-mail: fukui@somec.org; ORCID: 0009-0008-7122-522X).}%
\ifdefined\SUBMITTED\thanks{This work has been submitted to the IEEE for possible publication. Copyright may be transferred without notice, after which this version may no longer be accessible.}\fi}
\ifdefined\SUBMITTED
\markboth{Manuscript submitted to IEEE Transactions on Information Theory}{Fukui: Optimal reduced representatives of quadratic APN permutations}
\else
\markboth{Preprint, version 1.0.3; prepared for IEEE Transactions on Information Theory (not submitted)}{Fukui: Optimal reduced representatives of quadratic APN permutations}
\fi
\maketitle
\begin{abstract}
The formal derivative of the reduced polynomial representative of a map on a binary finite field depends on the
coordinates, although its differential uniformity does not. This matters for lifting permutations to Galois rings: a
polynomial over a Galois ring of characteristic a power of two greater than two permutes the ring exactly when its
reduction permutes the residue field and the derivative of the reduction has no zero there. We prove that for every odd $n\ge5$, every quadratic almost perfect
nonlinear permutation $Q$ of $\F_{2^n}$ admits invertible binary-linear maps $L$ and $A$ such that
the derivative of the reduced representative of $L\circ Q\circ A$ has no zero and each of its nonempty fibres has
exactly two elements. The range is exact: this fails for $n=3$, as computed earlier, and for even $n$ there are no such
permutations. For $n\ge9$, a linearization and a projectively injective cofactor map reduce the problem to counting
binary matrices: more than 5\% of the invertible inputs $A$ succeed, more than half for $n\ge15$, without any
classification. For $n=5,7$ we use two independent classifications and finite certificates. Consequently every such permutation is linearly equivalent to a counterexample to the reduced
critical-point conjecture of Bartoli and St\u{a}nic\u{a} that lifts to
permutation polynomials of optimal differential uniformity over every Galois ring with this residue field. For odd $n\ge9$ the existence proof is checked in the Lean~4 proof assistant.

\end{abstract}
\section{Introduction}\label{sec:intro}

\IEEEPARstart{A}{lmost} perfect nonlinear (APN) permutations of $\F_{2^n}$~\cite{Nyberg,Carlet} are studied through properties that do not
change under suitable changes of coordinates: the differential spectrum, the Walsh spectrum, and the equivalence
class under extended affine (EA) or CCZ equivalence. Polynomial lifting to Galois rings introduces a requirement of a
different kind. By a classical criterion (see~\cite{McDonald,Wan}; in this form~\cite[Theorem~2]{RS}), a polynomial over the Galois ring $\GR(2^k,n)$,
$k\ge2$, permutes the ring if and only if its reduction permutes $\F_{2^n}$ and the formal derivative of the reduction
has no zero on $\F_{2^n}$. The formal derivative of a representing polynomial is not an invariant of the map, and it
changes under linear changes of input and output even when the differential uniformity does not. This raises a
problem of choosing a representative.

\begin{quote}
\emph{Can every quadratic APN permutation in odd dimension be moved, by invertible $\F_2$-linear changes of input and
output, to a map whose reduced representative has a derivative without zeros and with the smallest possible
fibres?}
\end{quote}

We answer this affirmatively for every odd $n\ge5$ (Main Theorem below). The definitions are collected in
Section~\ref{sec:prelim}; we recall the essentials. Every map $f:F\to F$, $F=\F_{2^n}$, has a unique representing
polynomial $\red(f)\in F[X]$ of degree less than $2^n$, its \emph{reduced representative}, and $\red(f)'$ denotes its
formal derivative. A quadratic APN permutation $P$ is \emph{optimal} if $\red(P)'$ has no zero on $F$ and every
nonempty fibre of $\red(P)'$ has exactly two elements (Definition~\ref{def:opt}). For quadratic $P$ the map
$\red(P)'$ is affine; if it has no zero, it is not surjective, so its nonempty fibres have at least two elements.
Optimal maps are thus those whose derivative has no zero and the smallest fibres compatible with that.

\subsection{Main result}

\begin{mainthm}
Let $n\ge5$ be odd and let $Q$ be a quadratic APN permutation of $\F_{2^n}$. Then there are invertible $\F_2$-linear
maps $L,A$ of $\F_{2^n}$ such that $P=L\circ Q\circ A$ is optimal: $\red(P)'$ has no zero on $\F_{2^n}$, and there is
$r\ne0$ such that $\red(P)'(x)=\red(P)'(y)$ if and only if $y\in\{x,x+r\}$.
\end{mainthm}

The conclusion is stated for linear equivalence, which is finer than EA-equivalence; in particular the EA class
of every quadratic APN permutation in odd dimension $n\ge5$ contains an optimal representative. The quantifier ranges
over all quadratic APN permutations, including classes that have not been found; no list of families is used for
$n\ge9$.

The range of $n$ is exact. For even $n$ there is nothing to prove: an APN function of $\F_{2^n}$, $n$ even, all of
whose component functions are plateaued is not a permutation~\cite[Corollary~3(iii)]{BCCL}, and the components of a
quadratic function are quadratic, hence plateaued; so the class is empty for even $n$, the statement below holds
there vacuously, and nothing is constructed. For $n=1$ both permutations of $\F_2$ are quadratic in the sense of
Definition~\ref{def:quad} and optimal. For $n=3$ the statement is
false: no quadratic APN permutation of $\F_8$ has a reduced representative whose derivative is free of zeros
(Proposition~\ref{prop:nthree}; the same exhaustive result was reported earlier in~\cite{Stone-cx}). Thus every
quadratic APN permutation of $\F_{2^n}$ has an optimal linear equivalent if and only if $n\ne3$, vacuously for even $n$. Nothing is claimed for APN permutations that are not quadratic, or for quadratic APN functions
whose EA class contains no permutation.

Bartoli and St\u{a}nic\u{a}~\cite{BS-lift} showed that a reduced form of the lifting conjecture of~\cite{RS} is
equivalent to the following statement.

\begin{quote}
\emph{Reduced Critical-Point Conjecture} \cite[Conjecture~2]{BS-lift}. For every $n$, the reduced
representative of every APN permutation of $\F_{2^n}$ has a zero of its formal derivative in $\F_{2^n}$.
\end{quote}

\noindent For $n=1$ the identity map of $\F_2$ is a trivial counterexample, so the substantive cases are $n\ge3$.
Bartoli and St\u{a}nic\u{a} proved it for power maps whose reduced representative is $X^e$ with $e>1$, and for permutation
polynomials of degree three~\cite{BS-lift}. For such a power map the derivative vanishes at $0$, so the choice of
representative is essential for any counterexample.

\begin{corollary}\label{cor:bs}
For every odd $n\ge5$, every quadratic APN permutation of $\F_{2^n}$ is linearly equivalent to a counterexample to
the Reduced Critical-Point Conjecture.
\end{corollary}

\noindent Explicit counterexamples were posted independently by Stone~\cite{Stone-cx} after~\cite{BS-lift}
appeared: $\F_2$-affine input changes of APN power maps, for instance $(x+x^8+x^{16})^3$ on $\F_{2^5}$, whose
reduced representative has derivative $X^8+X^2+1$, with further examples for $n=7$ and $n=9$ given as data; no
general statement is claimed there. Corollary~\ref{cor:bs} concerns every quadratic APN permutation in every odd
dimension $n\ge5$, including classes that are not known. For $n=3$ the conjecture holds for quadratic APN
permutations, as already reported in~\cite{Stone-cx}; Proposition~\ref{prop:nthree} re-verifies this independently.

\begin{corollary}\label{cor:ring}
Let $P$ be optimal and $k\ge2$. Every polynomial over $\GR(2^k,n)$ whose coefficientwise reduction is $\red(P)$ is a
permutation polynomial of $\GR(2^k,n)$ with differential uniformity $\Delta^R=2\cdot2^{n(k-1)}$ (defined with the
difference of the ring, Section~\ref{sec:prelim}), the smallest value possible
for a permutation polynomial~\cite[Theorem~3]{RS}. Consequently, for every odd $n\ge5$, every quadratic APN
permutation of $\F_{2^n}$ is linearly equivalent to a map whose reduced representative lifts to such a permutation
polynomial of every $\GR(2^k,n)$.
\end{corollary}

\noindent Corollary~\ref{cor:ring} uses the permutation criterion and the differential estimates
of~\cite{RS} (Section~\ref{sec:cons}). R{\o}njom and Sandrib already give permutation polynomials with this
differential uniformity, from Dickson polynomials~\cite[Theorem~6]{RS}; their residue maps are not APN
(Proposition~\ref{prop:dickson}). The new point of Corollary~\ref{cor:ring} is that the residue map is an APN
permutation, in every class. Nothing is claimed about the security of any cipher.

\subsection{Method}

Fix $Q$ and an invertible $\F_2$-linear input $A$, and put $g=Q\circ A$. The outputs $L$ for which $\red(L\circ g)'$
is constant form a space $W_A$ over $F$; let $U_A$ be the common kernel of these maps, $s(A)=\dim_{\F_2}U_A$, and let
$\eta(A)\in\{0,1\}$ record whether the constant can be nonzero (Section~\ref{sec:prelim}).

\emph{The output at a fixed input} (Section~\ref{sec:layer1}). Theorem~\ref{thm:twoinv} shows that, at the input $A$,
three explicit constructions (two change the output, one composes the input with a shear) produce an optimal
$L\circ Q\circ A'$ whenever $s(A)=0$, or $s(A)=1$ and $\eta(A)=1$. This is a sufficient condition. The ingredients are a maximal-rank statement for left $F$-spaces of $\F_2$-linear maps (a known fact
on rank-metric codes~\cite{ABNR,PSZ}), the almost bent property of quadratic APN functions in odd
dimension~\cite{CCZ98}, and an elementary lemma on pairs of trace values.

\emph{The choice of the input} (Sections~\ref{sec:layer2}--\ref{sec:small}). For $n\ge9$ we show that a suitable input
always exists (Theorem~\ref{thm:universal}). Encode $A$ by the vector $\zeta_A\in F^n$. The matrix $N_A$, and
hence $W_A$ and $s(A)$, is determined by a matrix $G(\zeta_A)$ that is linear in $\zeta_A$, and the label $p_A$ is
an explicit polynomial expression in $\zeta_A$ (Theorem~\ref{thm:lin}). When a vector $\omega(\zeta_A)$ of
polynomials does not vanish, $W_A=F\omega(\zeta_A)$ and $s(A)=n-\rk\omega(\zeta_A)$, where $\rk$ is the rank of the
corresponding $\F_2$-linear map. The key step is that $[\zeta]\mapsto[\omega(\zeta)]$ is injective on projective classes outside
the zero set of an explicit nonzero polynomial $\Phi R$ of degree at most $n(n-1)$ (Lemma~\ref{lem:inj}). The
nonvanishing of $\Phi$ and $R$ is certified by evaluation at Moore points, using only the permutation property and
the almost bent property. Failing inputs are therefore counted by binary matrices of rank at most $n-2$ plus the
zeros of $\Phi R$. The former are fewer than the invertible matrices, and the Schwartz--Zippel bound controls the
latter for $n\ge9$. The argument is a pigeonhole comparison; it uses no character-sum estimate and no
classification. For $n=5,7$ the bound does not close, and we use the classifications of Brinkmann and
Leander~\cite{BrinkmannLeander}, Yoshiara~\cite{Yoshiara12,Yoshiara16} and Kalgin and Idrisova~\cite{KI}, together
with explicit input witnesses and $485$ inconsistency certificates (Section~\ref{sec:small}).

\subsection{Relation to previous work}

Table~\ref{tab:prior} compares the Main Theorem with the closest statements we know. The known ingredients of the
proof are the Moore determinant~\cite{Moore}, the almost bent property~\cite{CCZ98}, the Schwartz--Zippel
bound~\cite{Ore,DeMilloLipton,Zippel,Schwartz}, counts of binary matrices by rank, maximal-rank statements for
rank-metric codes~\cite{ABNR,PSZ}, the ring criterion and bounds of~\cite{RS}, and, for $n=5,7$, the
classifications cited above. What we add is the following.
\begin{itemize}
\item The reduction of the representative problem to the choice of one projective vector $[\zeta_A]$, through two
invariants of the input (Theorem~\ref{thm:twoinv}) and a linearization in $\zeta_A$ (Theorem~\ref{thm:lin}).
\item The projective injectivity of the cofactor map outside an explicit hypersurface (Lemma~\ref{lem:inj}), which
turns the existence of a good input into a count of binary matrices by rank.
\item A statement for every class, including classes that are not known, with a quantitative form: for odd $n\ge9$
at least a fraction $\pi_n$ of all invertible inputs works ($\pi_n$ is a lower bound), $\pi_n>0.054$ for every odd $n\ge9$ and $\pi_n>1/2$ for $n\ge15$
(Corollary~\ref{cor:proportion}); an optimal representative is therefore found by independent uniform sampling of
inputs (Remark~\ref{rem:algorithm}).
\item The exact range of dimensions, all $n\ne3$: the Main Theorem for odd $n\ge5$ together with~\cite{BCCL} (even
$n$, empty) and the exhaustive result for $n=3$ (Proposition~\ref{prop:nthree}, reported earlier in~\cite{Stone-cx}).
Of these, only the part $n\ge5$ is new.
\end{itemize}
Explicit examples, such as those of~\cite{Stone-cx} and of Appendix~\ref{app:gold}, cover particular classes; the
Main Theorem covers all of them. The contribution is thus not the existence of counterexamples
to~\cite[Conjecture~2]{BS-lift}, which~\cite{Stone-cx} already exhibits, but a statement for every class, including
classes that are not known, obtained by a method (projective injectivity of the cofactor map and a count of binary
matrices by rank) that uses no classification for $n\ge9$.

\begin{table}[t]\centering\small
\caption{Closest statements in the literature and the present result.}\label{tab:prior}
\begin{tabular}{@{}>{\raggedright}p{0.2\textwidth}>{\raggedright}p{0.25\textwidth}>{\raggedright\arraybackslash}p{0.45\textwidth}@{}}
\toprule
Source & Hypotheses & Conclusion \\
\midrule
\cite[Conj.~2]{BS-lift} & APN permutation, reduced representative & the derivative has a zero (conjectured) \\
\cite{BS-lift} & power map with reduced representative $X^e$, $e>1$; or permutation polynomial of degree $3$ & the derivative has a zero (proved) \\
\cite{Stone-cx} & $\F_2$-affine input changes of APN power maps; $n=5,7,9$ & explicit maps whose derivative has no zero (examples; no general statement); $n=3$: exhaustive check that the conjecture holds \\
\cite[Thms.~2,~3; Props.~4,~5]{RS} & polynomial over $\GR(2^k,n)$, $k\ge2$ & permutation criterion; bounds on differential uniformity \\
\cite[Thm.~6]{RS} & Dickson polynomial $g_5$, odd $n$ & optimal permutation polynomial; residue map not APN
(Prop.~\ref{prop:dickson}) \\
this paper & quadratic APN permutation, odd $n\ge5$ & some $L\circ Q\circ A$ is optimal; APN residue map with optimal lifts \\
\bottomrule
\end{tabular}
\end{table}

Appendix~\ref{app:gold} gives, independently of the Main Theorem, closed-form optimal representatives in every Gold
class; further explicit families and certificates are collected in Supplement~S4.

\subsection{Formal verification}

The Main Theorem for odd $n\ge9$ is formalized in the Lean~4 proof assistant with the Mathlib library and checked by
the Lean kernel, with no axioms beyond the three standard ones (Section~\ref{sec:lean}). The scope is as follows.
\begin{itemize}
\item Odd $n\ge9$: formalized; the formal statement concerns the actual objects (quadratic APN bijections of a finite
field of order $2^n$ and their Lagrange representatives).
\item $n=5,7$: not formalized; the proof rests on the cited classifications and finite certificates, and was
rechecked on a second, independently computed classification.
\item $n=3$ (Proposition~\ref{prop:nthree}; reported earlier in~\cite{Stone-cx}) and
Corollary~\ref{cor:proportion}: not formalized; exhaustive computation in two independent implementations, exact
rational arithmetic for finite $n$, and a rational interval for the limit.
\item Corollaries~\ref{cor:bs} and~\ref{cor:ring}: not formalized; the ring statements rest on~\cite{RS}.
\end{itemize}
The formal proof follows the proof given here; the paper is written to be read without it. The Lean development,
the certificates, the data and the scripts are publicly available (end of Section~\ref{sec:lean}).

\subsection{Organization}
Section~\ref{sec:prelim} gives definitions and notation. Section~\ref{sec:layer1} treats a fixed input.
Section~\ref{sec:layer2} proves the linearization in $\zeta_A$, Section~\ref{sec:cofactor} the cofactor dictionary
and the two nonvanishing polynomials, Section~\ref{sec:mainproof} the existence of a good input for $n\ge9$ and the
Main Theorem, and Section~\ref{sec:small} the dimensions $3$, $5$ and $7$. Section~\ref{sec:cons} proves the corollaries,
and Section~\ref{sec:lean} describes the formalization. Appendix~\ref{app:pfaff} gives a sharper degree bound and
Appendix~\ref{app:gold} the explicit Gold family. The supplementary material contains the correspondence between
statements and formal declarations, reproduction instructions, data and certificates (Supplement~S1--S4), and a
table of the identifiers used in the data archive (Supplement~S5).

\section{Definitions and notation}\label{sec:prelim}

\subsection{Maps of a finite field}
Throughout, $n\ge1$, $q=2^n$, $F=\F_q$, and $\Tr=\Tr_{F/\F_2}$. Every map $f:F\to F$ is represented by a unique
polynomial of degree less than $q$, the Lagrange interpolant; we call it the \emph{reduced representative} and write
$\red(f)\in F[X]$. Reduction is $F$-linear: $\red(\alpha f+\beta g)=\alpha\,\red(f)+\beta\,\red(g)$ for
$\alpha,\beta\in F$. The formal derivative of $\red(f)$ is denoted $\red(f)'$; the derivative of $X^e$ is $eX^{e-1}$,
which vanishes for even $e$. The \emph{binary weight} of an exponent is the number of ones in its binary expansion.

\begin{definition}[differential uniformity; APN]\label{def:apn}
For $a,b\in F$ let $\delta_f(a,b)=\#\{x\in F:f(x+a)+f(x)=b\}$. The map $f$ is \emph{almost perfect nonlinear}
(APN)~\cite{Nyberg} if $\delta_f(a,b)\le2$ for all $a\ne0$ and all $b$.
\end{definition}

\begin{definition}[quadratic; polar form]\label{def:quad}
The \emph{polar form} of $f$ is $B_f(x,y)=f(x+y)+f(x)+f(y)+f(0)$. The map $f$ is \emph{quadratic} (of algebraic
degree at most $2$) if $B_f$ is $\F_2$-bilinear; equivalently, every exponent of $\red(f)$ has binary weight at most
$2$.
\end{definition}

For quadratic $f$ and $a\ne0$, the map $x\mapsto f(x+a)+f(x)=B_f(x,a)+f(a)+f(0)$ is affine, so $f$ is APN iff the
$\F_2$-linear map $B_f(\cdot,a)$ has kernel $\{0,a\}$ for every $a\ne0$. All maps below are considered as maps
$F\to F$; ``linear'' means $\F_2$-linear, and $\GL(F)$ is the group of invertible $\F_2$-linear maps of $F$.

\begin{definition}[equivalences]\label{def:equiv}
Maps $f,g:F\to F$ are \emph{linearly equivalent} if $g=L\circ f\circ A$ with $L,A\in\GL(F)$, and
\emph{EA-equivalent} (extended affine) if $g=L\circ f\circ A+C$ with $L,A$ affine bijections and $C$ affine.
CCZ-equivalence is equivalence of graphs under an affine bijection of $F\times F$~\cite{CCZ98,BCP}.
\end{definition}

The \emph{Gold functions} are the power maps $x^{2^k+1}$ with $\gcd(k,n)=1$; they are quadratic APN functions, and
permutations when $n$ is odd. Linear equivalence implies EA-equivalence, which implies CCZ-equivalence. Linear equivalence preserves being
quadratic, APN and bijective. The formal derivative of the reduced representative is not preserved.

\begin{definition}[optimal]\label{def:opt}
A quadratic APN permutation $P$ of $F$ is \emph{optimal} if $\red(P)'$ has no zero on $F$ and every nonempty fibre of
the map $x\mapsto\red(P)'(x)$ on $F$ has exactly two elements.
\end{definition}

\begin{remark}[optimal maps]\label{rem:optaffine}
Let $P$ be quadratic. The exponents of $\red(P)$ have weight at most $2$, and the odd ones are $1$ and $1+2^j$,
$1\le j\le n-1$; so $\red(P)'(X)=c+K(X)$ with $c\in F$ and $K(X)=\sum_{j=1}^{n-1}k_jX^{2^j}$ additive. The fibres of
$c+K$ are the cosets of $\ker K$. Hence $\red(P)'$ has no zero and all its nonempty fibres have exactly two elements
iff $\rank_{\F_2}K=n-1$ and $c\notin\im K$; then $\ker K=\{0,r\}$ with $r\ne0$, and $\red(P)'(x)=\red(P)'(y)$ iff
$y\in\{x,x+r\}$. In words: $\red(P)'$ is exactly two-to-one onto its image, which is a hyperplane coset not
containing $0$. It is not a map of $F$ onto $F$. If moreover $P$ is an APN permutation, this derivative condition is
exactly optimality. The derivative condition does not imply the APN property: $X^5+X^3+X$ on $\F_{2^5}$ is a
quadratic permutation whose derivative $1+X^2+X^4$ satisfies it, but it is not APN
(Proposition~\ref{prop:dickson} with $u=1$; its differential uniformity is $4$).
\end{remark}

\subsection{Galois rings}
For $k\ge2$, $\GR(2^k,n)$ is the Galois ring of characteristic $2^k$ with residue field $F$. A polynomial over
$\GR(2^k,n)$ \emph{lifts} $f\in F[X]$ if its coefficientwise reduction modulo $2$ is $f$. The differential uniformity
of a map $\widehat f$ of the ring is defined with the difference of the ring, as in~\cite[Section~3]{RS}:
\[
 \delta^R_{\widehat f}(a,b)=\#\{x\in\GR(2^k,n):\widehat f(x+a)-\widehat f(x)=b\},\qquad
 \Delta^R_{\widehat f}=\max_{a\ne0,\ b}\delta^R_{\widehat f}(a,b).
\]
On $F$ the difference and the sum coincide, because $F$ has characteristic $2$, so the same formula on $F$
gives $\delta_f$ of Definition~\ref{def:apn}; on $\GR(2^k,n)$ with $k\ge2$ they differ (for the identity map of $\GR(4,n)$ the largest
set $\{x:(x+a)+x=b\}$ has $2^n$ elements, while $\{x:(x+a)-x=a\}$ is the whole ring).

\subsection{Notation for the proof}
From now on $n\ge5$ is odd and $Q$ is a quadratic APN permutation of $F$.

\emph{Outputs.} For $\ell\in F^n$ put $L_\ell(y)=\sum_{i=0}^{n-1}\ell_iy^{2^i}$ and $m(y)=(y^{2^i})_{i=0}^{n-1}$, so
that $L_\ell(y)=m(y)\cdot\ell$. The map $\ell\mapsto L_\ell$ is a bijection from $F^n$ onto the $\F_2$-linear maps
$F\to F$. Write $\rk\ell=\dim_{\F_2}L_\ell(F)$.

\emph{Interpolation data.} For a quadratic map $g$ let $N=N_g$ be the $(n-1)\times n$ matrix with entries
$N_{j,i}=[X^{1+2^j}]\red(g^{2^i})$, $1\le j\le n-1$, $0\le i<n$, and $p=p_g\in F^n$ the vector with
$p_i=[X^1]\red(g^{2^i})$. Then
\[\red(L_\ell\circ g)'=p\cdot\ell+K_\ell,\qquad K_\ell(X)=\sum_{j=1}^{n-1}(N\ell)_jX^{2^j}.\]
Put $W=\ker_FN$ (a left $F$-space of outputs, via $\ell\mapsto L_\ell$), $d=\dim_FW$,
$U=\{y\in F:L_\ell(y)=0\text{ for all }\ell\in W\}$, $s=\dim_{\F_2}U$, and $\eta=1$ if $p\notin\row_FN$, $\eta=0$
otherwise (equivalently, $\eta=1$ iff $\ell\mapsto p\cdot\ell$ is not identically zero on $W$). For an input
$A\in\GL(F)$ we write $N_A,p_A,W_A,U_A,d(A),s(A),\eta(A)$ for the data of $g=Q\circ A$. Thus $W_A$ is the space of
outputs $L$ for which $\red(L\circ Q\circ A)'$ is constant, $s(A)$ the dimension of their common kernel, and
$\eta(A)=1$ iff that constant can be nonzero.

\emph{Inputs.} For $A(x)=\sum_k\alpha_kx^{2^k}$ let $\zeta_A\in F^n$, $(\zeta_A)_c=\alpha_{-c}^{2^c}$ (indices modulo
$n$; this is column $0$ of the Dickson matrix $\mathcal D_A$ of $A$, Section~\ref{sec:layer2}). The map $A\mapsto\zeta_A$ is a bijection from the
$\F_2$-linear maps $F\to F$ onto $F^n$.

\emph{The polar tensor of $Q$.} Let $c_e$ be the reduced coefficients of $Q$. Put $\Gamma_{a,c,i}=c_{2^{a-i}+2^{c-i}}^{\,2^i}$
for $a\ne c$ and $\Gamma_{a,a,i}=0$ (indices modulo $n$), and $\Gamma(\xi,\chi)_i=\sum_{a,c}\xi_a\chi_c\Gamma_{a,c,i}$ for
$\xi,\chi\in F^n$. Let $G(\zeta)$ be the $n\times n$ matrix with rows $\Gamma(\zeta,e_c)$, and $\bG_\ell$ the
alternating matrix $\bG_\ell[a,c]=\sum_i\Gamma_{a,c,i}\ell_i$; then $G(\zeta)\ell=\bG_\ell\zeta$. The vector
$\omega(Z)$ of polynomials and the polynomial $\Phi(Z)$ are defined in Theorem~\ref{thm:lin}.

We distinguish two ranks: $\rank$ is the rank of a matrix over $F$, and $\rk\ell$ is the $\F_2$-rank of the map
$L_\ell$.

\section{The output at a fixed input}\label{sec:layer1}

Throughout this section $n\ge5$ is odd unless stated otherwise. We work at a fixed input: $Q$ is a quadratic APN
permutation, $A$ an invertible $\F_2$-linear map and $g=Q\circ A$. The section gives a sufficient condition, in
terms of two invariants of $g$, for an optimal representative obtained from $g$ by a change of output and, in the
zero-derivative branch, one additional input shear $A\mapsto A\circ A_t$ (Theorem~\ref{thm:twoinv}). The condition is
not claimed to be necessary for the existence of an optimal output at the input $A$. Several ingredients are
standard facts on rank-metric codes and almost bent functions; we give short proofs so that the argument is
self-contained, and we indicate the known sources. No character-sum estimate is used in this section.

\subsection{The interpolation data}\label{ss:data}
Recall from Section~\ref{sec:prelim} that $\red(L_\ell\circ g)'=p\cdot\ell+K_\ell$ with
$K_\ell(X)=\sum_{j=1}^{n-1}(N\ell)_jX^{2^j}$. We recall why both parts are $F$-linear in $\ell$, with
no Frobenius twist. Since $L_\ell\circ g=\sum_i\ell_i\,g^{2^i}$ and reduction is $F$-linear, every coefficient of
$\red(L_\ell\circ g)$ has the form $\sum_i\ell_ir_i$ with $r_i\in F$ independent of $\ell$. The exponents of
$\red(g^{2^i})$ have binary weight $\le2$; the formal derivative of $X^e$ is $eX^{e-1}$, nonzero only for odd $e$,
and among weights $\le2$ these are $e=1$ (giving the constant) and $e=1+2^j$, $1\le j\le n-1$ (giving the
coefficient of $X^{2^j}$). Hence $W=\ker_FN$ is an $F$-subspace of $F^n$, and since $\ell\mapsto L_\ell$ is
$F$-linear for the \emph{left} action $(\alpha L)(y)=\alpha L(y)$, $\{L_\ell:\ell\in W\}$ is a left $F$-subspace of
$\mathrm{Hom}_{\F_2}(F,F)$. We identify $\ell$ with $L_\ell$ throughout. \emph{CD} abbreviates \emph{constant derivative}: \emph{CD at $A$} means that some $\ell\in W_A$
has $L_\ell$ bijective; a \emph{CD output} is any $\ell\in W_A$, with no normalization of the constant; a
\emph{common root} is a common \emph{nonzero} root in $F$.

\subsection{Translation invariance}\label{ss:transl}
\begin{lemma}\label{lem:L1red}
For every $f:F\to F$ and $b\in F$, $\red(f(\cdot+b))(X)=\red(f)(X+b)$; hence $\red(f(\cdot+b))'(X)=\red(f)'(X+b)$.
In particular, if $\red(f)'$ is a constant $\delta$ (zero or not), then $\red(f(\cdot+b))'=\delta$ for every $b$.
\end{lemma}
\begin{proof}
$\red(f)(X+b)$ has degree $<q$ and takes the value $f(x+b)$ at every $x$; uniqueness of the representative of
degree $<q$ gives the first identity. The formal derivative commutes with the substitution $X\mapsto X+b$.
\end{proof}

\begin{proposition}\label{prop:L1Np}
Let $g$ be quadratic and $g_b=g(\cdot+b)$. Then
\[N_{g_b}=N_g,\qquad p_{g_b}-p_g=N_g^{\top}\beta\in\row_FN_g,\qquad \beta=(b^{2^j})_{j=1}^{n-1}.\]
Consequently $\ell\mapsto p\cdot\ell$ is unchanged on $W=\ker_FN_g$; in particular $\eta$ is unchanged.
\end{proposition}
\begin{proof}
By Lemma~\ref{lem:L1red} applied to $f=L_\ell\circ g$,
$\red(L_\ell\circ g_b)'(X)=p\cdot\ell+K_\ell(X+b)=\bigl(p\cdot\ell+K_\ell(b)\bigr)+K_\ell(X)$, since $K_\ell$ is additive.
The non-constant part is unchanged, so $N$ is unchanged; the constant changes by
$K_\ell(b)=\sum_j(N\ell)_jb^{2^j}=(N^\top\beta)\cdot\ell$, an $F$-linear functional of $\ell$ lying in $\row_FN$
and vanishing on $\ker_FN$.
\end{proof}

\subsection{Maximal rank in a left \texorpdfstring{$F$}{F}-space}\label{ss:maxrank}
In \S\ref{ss:maxrank}--\ref{ss:label}, $W\subseteq\mathrm{Hom}_{\F_2}(F,F)$ is any left $F$-subspace, $d=\dim_FW$,
$U=\bigcap_{L\in W}\ker L$ and $s=\dim_{\F_2}U$ ($U=F$ if $W=0$). APN is not used. For $y\in F$ put
$\mathrm{ev}_y:W\to F$, $L\mapsto L(y)$.

\begin{lemma}\label{lem:L1ev}
(a) $\mathrm{ev}_y$ is $F$-linear. (b) $\mathrm{ev}_y=0$ iff $y\in U$. (c) $\mathrm{ev}_{y+u}=\mathrm{ev}_y$ for $u\in U$.
Hence $y\mapsto\mathrm{ev}_y$ factors through $F/U$, and the nonzero functionals of the form $\mathrm{ev}_y$ are
indexed by at most $2^{n-s}-1$ nonzero cosets.
\end{lemma}
\begin{proof}
(a) $\mathrm{ev}_y(\alpha L+\beta L')=\alpha L(y)+\beta L'(y)$ by the definition of the left action; this, together
with $|W|=q^d$ and $|\ker\mathrm{ev}_y|=q^{d-1}$ in Theorem~\ref{thm:L1max}, is where $F$-closure enters, and it is
essential (Example~\ref{ex:L1neg}). (b) is the definition of $U$. (c) $L(y+u)=L(y)+L(u)=L(y)$ for $L\in W$.
\end{proof}

\begin{theorem}[maximal rank]\label{thm:L1max}
$\max_{L\in W}\rk L=n-s$. More precisely, if $W\ne0$, the set $\mathcal G=\{L\in W:\ker L=U\}$ satisfies
$|\mathcal G|\ge q^d-(2^{n-s}-1)q^{d-1}\ge q^{d-1}>0$.
\end{theorem}
\begin{proof}
Every $L\in W$ kills $U$, so $\rk L\le n-s$. If $W=0$ then $s=n$ and the maximum is $0=n-s$. Let $W\ne0$. Then
$L\notin\mathcal G$ iff $L(y)=0$ for some $y\notin U$, i.e.\ $L\in\ker\mathrm{ev}_y$. By Lemma~\ref{lem:L1ev} each such
$\ker\mathrm{ev}_y$ is an $F$-hyperplane of $W$ ($q^{d-1}$ elements), and there are at most $2^{n-s}-1\le q-1$
distinct ones. Their union has at most $(2^{n-s}-1)q^{d-1}<q^d$ elements. Any $L\in\mathcal G$ has kernel exactly $U$,
so rank $n-s$.
\end{proof}

\begin{corollary}\label{cor:L1bij}
(a) If $W=F\cdot L_1$ with $L_1\ne0$, then every nonzero element of $W$ has kernel exactly $U=\ker L_1$.
(b) $W$ contains a bijection iff $s=0$, i.e.\ iff the maps in $W$ have no common nonzero root. Applied to
$W_A=\ker_FN_A$: CD holds at $A$ iff $s(A)=0$ and, whether or not CD holds, the
largest rank of a CD output at $A$ is exactly $n-s(A)$.
\end{corollary}
\begin{proof}
(a) $\ker(\alpha L_1)=\ker L_1$ for $\alpha\ne0$. (b) $s=0$ gives rank $n$ by Theorem~\ref{thm:L1max}; conversely a
bijection has trivial kernel, so $U=0$. $W_A$ is a left $F$-subspace by \S\ref{ss:data}.
\end{proof}

\begin{example}[$F$-closure cannot be dropped]\label{ex:L1neg}
Let $n\ge2$, $u\in F^\times$, and $W'=\{y\mapsto\Tr(ay)\,u:a\in F\}$. This $\F_2$-space has common kernel $0$, but
every element has rank $\le1$. It is not left $F$-closed, and the proof breaks at Lemma~\ref{lem:L1ev}(a):
$\ker\mathrm{ev}_y$ is then an $\F_2$-hyperplane, and $2^n-1$ of them cover $W'$.
\end{example}

\begin{remark}
The argument is the elementary fact that at most $q-1$ proper $F$-hyperplanes cannot cover an $F$-space of
dimension $\ge1$ (cf.\ \cite{Clark}), plus the factorization through $F/U$; no novelty is claimed. In coding
language the maximal-rank statement is \cite[Proposition~3.11]{ABNR}, and a count of maximum-rank elements is
in \cite[Theorems~3.2 and~3.9]{PSZ}.
\end{remark}

\subsection{Maximal rank on a prescribed nonzero label}\label{ss:label}
\begin{theorem}[prescribed nonzero label]\label{thm:L1lab}
Let $W\ne0$, $\mathcal G$ as in Theorem~\ref{thm:L1max}, $c:W\to F$ a nonzero $F$-linear functional, $\delta\in F^\times$
and $H=\{L\in W:c(L)=\delta\}$. Then
\[|H\cap\mathcal G|\ \ge\ \begin{cases}1,&d=1,\\ q^{d-2}\bigl(q-q/2^s+1\bigr),&d\ge2;\end{cases}\]
in particular some $L\in W$ has $c(L)=\delta$ and $\ker L=U$ (rank $n-s$).
\end{theorem}
\begin{proof}
$H$ is an affine hyperplane of $W$, $|H|=q^{d-1}$, and $0\notin H$ since $\delta\ne0$.
\emph{$d=1$.} $W=FL_1$ and $c(\alpha L_1)=\alpha c(L_1)$ with $c(L_1)\ne0$, so $H=\{(\delta/c(L_1))L_1\}$, a nonzero
multiple of $L_1$, whose kernel is $U$.
\emph{$d\ge2$.} $L\in H\setminus\mathcal G$ iff $L\in H$ and $\mathrm{ev}_y(L)=0$ for some $y\notin U$. Fix $y\notin U$;
$\mathrm{ev}_y$ is a nonzero $F$-linear functional (Lemma~\ref{lem:L1ev}). (i) If $\mathrm{ev}_y=\lambda c$, then
$\lambda\ne0$ and $\mathrm{ev}_y\equiv\lambda\delta\ne0$ on $H$: no point of $H$ is excluded. (ii) If
$\mathrm{ev}_y$ and $c$ are $F$-linearly independent, $(c,\mathrm{ev}_y):W\to F^2$ is surjective and
$H\cap\ker\mathrm{ev}_y$ is the fibre over $(\delta,0)$, of size $q^{d-2}$. Since $\mathrm{ev}_y$ depends only on
$y+U$, at most $q/2^s-1$ distinct conditions occur. Hence $|H\setminus\mathcal G|\le(q/2^s-1)q^{d-2}$ and
$|H\cap\mathcal G|\ge q^{d-1}-(q/2^s-1)q^{d-2}=q^{d-2}(q-q/2^s+1)\ge q^{d-2}>0$.
\end{proof}

\begin{example}[$\delta=0$ is different]\label{ex:L1zero}
Let $n\ge2$, $W=F\,\mathrm{id}\oplus F\,\Tr$ ($d=2$), $c(a\,\mathrm{id}+b\Tr)=a$. Then $s=0$, but $\ker c=F\Tr$
consists of maps of rank $\le1$. For $\delta=0$ the relevant statement is Theorem~\ref{thm:L1max} for the
subspace $\ker c$ and its own common kernel, which can be larger than $U$.
\end{example}

\begin{corollary}[normalization for the kernel-line route]\label{cor:L1norm}
Let $c_A(\ell)=p_A\cdot\ell$. The hypothesis of Proposition~\ref{prop:outKL} at $A$ is: some $w\in W_A$ has
$c_A(w)\ne0$ and $\rk L_w\ge n-1$. If $c_A\not\equiv0$ on $W_A$ and $s(A)\le1$, this hypothesis holds; moreover,
for every $\delta\ne0$ the witness can be chosen with $c_A(w)=\delta$ and $\ker L_w=U_A$.
\end{corollary}
\begin{proof}
$W_A\ne0$ and $c_A|_{W_A}$ is a nonzero $F$-linear functional; Theorem~\ref{thm:L1lab} gives $w$ with
$c_A(w)=\delta$ and $\ker L_w=U_A$, so $\rk L_w=n-s(A)\ge n-1$.
\end{proof}

\subsection{The polar identity and minimum rank three}\label{ss:min3}
For $f\in F$ put $\ell_f=(f^{2^i})_i=m(f)$, so that $L_{\ell_f}=\Tr(f\,\cdot)$, and $K_f:=K_{\ell_f}$,
$c_f:=p\cdot\ell_f$. In this subsection and the next, $g$ is quadratic (every exponent of $\red(g)$ has binary
weight $\le2$) and need not be a permutation.

\begin{lemma}[polar identity]\label{lem:L1polar}
Let $n\ge3$, $\varphi=\Tr(fg)$ and $B_f(x,y)=\varphi(x+y)+\varphi(x)+\varphi(y)+\varphi(0)$. Then
$B_f(x,y)=\Tr(x\,K_f(y))$ for all $x,y$. Consequently $K_f$ is self-adjoint for $\langle x,y\rangle=\Tr(xy)$,
$\ker K_f=\operatorname{rad}B_f$, and $\rk K_f=n-\dim\operatorname{rad}B_f$.
\end{lemma}
\begin{proof}
$\varphi$ is $\F_2$-valued, so $\red(\varphi)^2$ reduces to $\red(\varphi)$; writing $\red(\varphi)=\sum_ec_eX^e$
this says $c_{2e\bmod(q-1)}=c_e^2$ for $1\le e\le q-2$ (for $n\ge3$ no weight-two exponent equals $q-1$). The map
$\varphi=\sum_i(fg)^{2^i}$ is again quadratic. A monomial of weight one contributes nothing to $B_f$; a monomial
$X^{2^a+2^b}$, $a\ne b$, contributes $x^{2^a}y^{2^b}+x^{2^b}y^{2^a}$. Hence
$B_f(x,y)=\sum_{a\ne b}c_{2^a+2^b}x^{2^a}y^{2^b}$ (ordered pairs). As in \S\ref{ss:data},
$K_f(y)=\sum_{j=1}^{n-1}c_{1+2^j}y^{2^j}$, and
\[\Tr(xK_f(y))=\sum_{k=0}^{n-1}\sum_{j=1}^{n-1}c_{1+2^j}^{2^k}x^{2^k}y^{2^{j+k}}
=\sum_k\sum_jc_{2^k+2^{j+k}}x^{2^k}y^{2^{j+k}},\]
using $c_{2^ke}=c_e^{2^k}$. The map $(k,j)\mapsto(a,b)=(k,k+j\bmod n)$ is a bijection onto the ordered pairs with
$a\ne b$, which gives $B_f$. Symmetry of $B_f$ gives self-adjointness, and nondegeneracy of the trace form gives
$\ker K_f=\{y:B_f(x,y)=0\ \forall x\}=\operatorname{rad}B_f$.
\end{proof}

\begin{lemma}[APN, $n$ odd]\label{lem:L1apn}
If $g$ is quadratic APN and $n$ is odd, then $\dim\operatorname{rad}B_f=1$, i.e.\ $\rk K_f=n-1$, for every $f\ne0$.
Moreover, for every $a\ne0$ exactly one $f\ne0$ has $a\in\ker K_f$; writing $\ker K_f=\{0,r_f\}$, the map
$f\mapsto r_f$ is a bijection $F^\times\to F^\times$.
\end{lemma}
\begin{proof}
For $a\ne0$ the map $\beta_a(x)=g(x+a)+g(x)+g(a)+g(0)$ is $\F_2$-linear. APN makes $x\mapsto g(x+a)+g(x)$
two-to-one, so $\ker\beta_a=\{0,a\}$ and $\operatorname{im}\beta_a$ is a hyperplane. Hence exactly one $f\ne0$ has
$\Tr(f\beta_a(x))=0$ for all $x$. Since $\Tr(f\beta_a(x))=B_f(x,a)$, for each $a\ne0$ exactly one $f\ne0$ has
$a\in\operatorname{rad}B_f$. So there are exactly $q-1$ pairs $(a,f)$ with $a,f\ne0$ and $a\in\operatorname{rad}B_f$.
Each $B_f$ is alternating, hence of even rank, and for odd $n$ its radical has dimension $\ge1$. Thus each of the
$q-1$ values of $f$ contributes at least one pair, and equality forces exactly one nonzero radical vector. The
last sentence restates the count.
\end{proof}
A function is \emph{almost bent} if all its nonzero component functions $\Tr(fg)$ have Walsh spectrum in
$\{0,\pm2^{(n+1)/2}\}$; for quadratic functions this is equivalent to $\dim\operatorname{rad}B_f=1$ for all $f\ne0$, which is
the statement above. So Lemma~\ref{lem:L1apn} is the classical fact that quadratic APN functions in odd dimension are almost
bent \cite{CCZ98}; the double
count is included for completeness.

\begin{theorem}[minimum rank three]\label{thm:L1min3}
Let $n\ge3$ be odd and $g$ quadratic APN (not necessarily a permutation). Every nonzero $L\in W_g$ has
$\rk L\ge3$.
\end{theorem}
\begin{proof}
An $\F_2$-linear $L$ of rank $r$ can be written $L(y)=\sum_{t=1}^ru_t\Tr(f_ty)$ with $u_1,\dots,u_r$ and
$f_1,\dots,f_r$ each $\F_2$-independent. As a left $F$-combination $L=\sum_tu_tL_{\ell_{f_t}}$, so its linear part is
$\sum_tu_tK_{f_t}$, and $L\in W_g$ iff this is $0$.
\emph{Rank one.} $u_1K_{f_1}=0$ with $u_1\ne0$ forces $K_{f_1}=0$, contradicting $\rk K_{f_1}=n-1\ge2$
(Lemma~\ref{lem:L1apn}).
\emph{Rank two.} $u_1K_{f_1}=u_2K_{f_2}$. Put $\alpha=u_1/u_2\notin\{0,1\}$; then $K_{f_2}=m_\alpha\circ K_{f_1}$ with
$m_\alpha(z)=\alpha z$. Taking adjoints ($m_\alpha^*=m_\alpha$, $K_{f_i}^*=K_{f_i}$ by Lemma~\ref{lem:L1polar}),
$K_{f_2}=K_{f_1}\circ m_\alpha$. Hence $m_\alpha$ commutes with $K_{f_1}$ and maps $\ker K_{f_1}=\{0,r\}$ into itself;
$\alpha r\ne0$ gives $\alpha r=r$ and $\alpha=1$, a contradiction.
\end{proof}

\begin{remark}
APN is essential: for $g=x^3+c\,\Tr(x^3)$ with $\Tr c=1$ the component $\Tr(g)$ vanishes, so $K_1=0$ and the rank-one
map $\Tr$ lies in $W_g$. Rank three occurs: for $g=x^3$ over $\F_{32}$ the $32768$ elements of $W_g$ have ranks
$\{0\!:\!1,\ 3\!:\!4805,\ 4\!:\!17298,\ 5\!:\!10664\}$ (direct computation; Supplement~S3).
\end{remark}

\subsection{A dimension constraint}\label{ss:dim}
\begin{theorem}\label{thm:L1dim}
Let $n\ge3$ be odd and $g$ quadratic APN. Then $W_g\ne0$, $\rank N\ge2$, and
\[s\le n-3,\qquad d+s\le n-2,\qquad s\le\rank N-2,\qquad \max_{L\in W_g}\rk L=n-s\ge d+2.\]
\end{theorem}
\begin{proof}
$N$ is $(n-1)\times n$, so $d=n-\rank N\ge1$ and $W\ne0$. A nonzero $L\in W$ kills $U$, so $3\le\rk L\le n-s$
(Theorem~\ref{thm:L1min3}); hence $s\le n-3$ and $m:=n-s-2\ge1$. Choose $\F_2$-independent $y_1,\dots,y_m$ in a
complement of $U$ and let $E:W\to F^m$, $L\mapsto(L(y_1),\dots,L(y_m))$, which is $F$-linear for the left action. If
$E(L)=0$ then $\ker L\supseteq U\oplus\langle y_1,\dots,y_m\rangle$, of dimension $n-2$, so $\rk L\le2$ and $L=0$.
So $E$ is injective and $d\le n-s-2$. With $d=n-\rank N$ this is $s\le\rank N-2$, hence $\rank N\ge2$; with
Theorem~\ref{thm:L1max} the maximal rank is $n-s\ge d+2$.
\end{proof}

\begin{remark}
Identify $L\in W$ with its evaluation vector $(L(b_1),\dots,L(b_n))$ for an $\F_2$-basis $(b_i)$ of $F$; this is an
$F$-linear injection, and the rank of the vector (over $\F_2$) equals $\rk L$. (The coefficient vector $\ell$ does
not have this property.) So $W$ is an $F$-linear rank-metric code of length $n$ and minimum distance $\ge3$, and
the rank-metric Singleton bound (see \cite{Sheekey}) gives $d\le n-2$, the case $s=0$. The bound $d+s\le n-2$ is
the same argument after factoring through $F/U$. No novelty is claimed.
\end{remark}

\subsection{Three output constructions}\label{ss:outputs}
The fixed-input criterion (Theorem~\ref{thm:twoinv}(iii)) ends with one of three explicit constructions, which we
prove here. They use only the following elementary lemma on trace values.

\begin{lemma}[four trace values]\label{lem:trace4}
Let $n\ge4$, $\gamma\in F^\times$ and $\beta,\varepsilon\in\F_2$. There is $z\in F^\times$ with $\Tr(z)=\beta$ and
$\Tr(\gamma/z)=\varepsilon$.
\end{lemma}
\begin{proof}
Put $M=q/2$, $T(X)=\sum_{i<n}X^{2^i}$, $\mathcal R(X)=\sum_{i<n}\gamma^{2^i}X^{M-2^i}$ and $\mathcal P(X)=X^M+\mathcal R(X)$. For $z\ne0$,
$\mathcal R(z)=z^M\Tr(\gamma/z)$ and $\mathcal P(z)=z^M(1+\Tr(\gamma/z))$. The exponents $M-2^i$ are distinct; $\mathcal R$ has degree $M-1$
with coefficients $\gamma$ at $X^{M-1}$ and $\gamma^2$ at $X^{M-2}$. $T$ is separable ($T'=1$), so the $M$ trace-one
elements are the roots of $H_1=T+1$ (monic of degree $M$, exponents $0$ and powers of two), and the $M-1$ nonzero
trace-zero elements are the roots of $H_0=T/X=\sum_iX^{2^i-1}$ (monic of degree $M-1$, exponents $0$ and odd numbers).
Suppose no $z$ exists.
\emph{$(\beta,\varepsilon)=(1,1)$:} $\mathcal R$ vanishes at $M$ points and has degree $<M$, so $\mathcal R=0$; but $[X^{M-1}]\mathcal R=\gamma\ne0$.
\emph{$(0,1)$:} $\mathcal R$ vanishes at the roots of $H_0$; $\mathcal R-\gamma H_0$ has degree $\le M-2$ and $M-1$ roots, so
$\mathcal R=\gamma H_0$; but $M-2\ge2$ is even, so $[X^{M-2}]H_0=0\ne\gamma^2$.
\emph{$(1,0)$:} $\mathcal P$ vanishes at the roots of $H_1$ and both are monic of degree $M$, so $\mathcal P=H_1$; but $M-1\ge3$ is odd,
so $[X^{M-1}]H_1=0\ne\gamma$.
\emph{$(0,0)$:} $H_0\mid\mathcal P$, so $\mathcal P=(X+u)H_0$, whose exponents are $2^i$ or $2^i-1$. For $n\ge4$, $M-2=2(2^{n-2}-1)$ is
even and not a power of two, so $[X^{M-2}]\mathcal P=0$; but $[X^{M-2}]\mathcal P=\gamma^2\ne0$.
\end{proof}
For $n=3$, $\gamma=1$, $(\beta,\varepsilon)=(0,0)$ the statement fails; we only use $n\ge5$.

\begin{lemma}[component derivatives]\label{lem:comp}
Let $g$ be quadratic APN, $n$ odd, $f\ne0$, and write $\red(\Tr(fg))'=c_f+K_f$ with $K_f$ additive
(\S\ref{ss:min3}). Then $K_f$ is self-adjoint for $\langle x,y\rangle=\Tr(xy)$, $\ker K_f=\{0,r_f\}$ with
$r_f\ne0$, and $\im K_f=\{v:\Tr(r_fv)=0\}$.
\end{lemma}
\begin{proof}
Self-adjointness and $\rk K_f=n-1$ are Lemmas~\ref{lem:L1polar} and~\ref{lem:L1apn}. For self-adjoint $K$,
$\im K=(\ker K)^\perp$ for the nondegenerate trace form.
\end{proof}

\begin{proposition}[nonzero constant]\label{prop:outC}
Let $n\ge5$ be odd and $h$ a quadratic APN permutation with $\red(h)'=\delta\in F^\times$. For every $\theta\in F^\times$
there is $\sigma\in F^\times$ such that $L(y)=\sigma y+\Tr(\theta y)$ is bijective and $L\circ h$ is optimal.
\end{proposition}
\begin{proof}
If $y\ne0$ and $\sigma y=\Tr(\theta y)$, then $\sigma y=1$ and $\Tr(\theta/\sigma)=1$; so $L$ is bijective iff $\Tr(\theta/\sigma)=0$. Reduction
is $F$-linear, so $\red(L\circ h)=\sigma\,\red(h)+\red(\Tr(\theta h))$ and $\red(L\circ h)'=(\sigma\delta+c_\theta)+K_\theta$. By
Lemma~\ref{lem:comp} the linear part has rank $n-1$, and the constant is outside its image iff
$\Tr(r_\theta\delta\sigma)=1+\Tr(r_\theta c_\theta)$. Lemma~\ref{lem:trace4} with $z=r_\theta\delta\sigma$, $\gamma=r_\theta\delta\theta\ne0$,
$\beta=1+\Tr(r_\theta c_\theta)$, $\varepsilon=0$ gives $\sigma\ne0$ satisfying both conditions (note $\Tr(\theta/\sigma)=\Tr(\gamma/z)$).
Then $L\circ h$ is a quadratic APN permutation, and its reduced derivative is $c+K$ with $\rk K=n-1$ and
$c\notin\im K$, i.e.\ root-free and exactly two-to-one.
\end{proof}

\begin{proposition}[zero derivative, input shear]\label{prop:outNF}
Let $n\ge5$ be odd and $S$ a quadratic APN permutation with $\red(S)'=0$. There is $t\ne0$ with $\Tr t=0$ such that
$S\circ A_t$ is optimal, $A_t(x)=x+\Tr(tx)$ (a bijection since $\Tr t=0$).
\end{proposition}
\begin{proof}
Put $H(x)=B_S(1,x)$ and $s_S=S(1)+S(0)$. By APN, $\ker H=\{0,1\}$, so $\im H=\{y:\Tr(\pi y)=0\}$ for some $\pi\ne0$;
and $s_S\notin\im H$, since $s_S=H(x)$ would give $S(x+1)=S(x)$. So $\Tr(\pi s_S)=1$. Since $\red(S)'=0$, every
exponent of $\red(S)$ is even, hence $H=\sum_{k=1}^{n-1}h_kX^{2^k}$ has no $X$-term; put $u=h_{n-1}$. For
$\ell=\Tr(tx)\in\F_2$, $S(x+\ell)=S(x)+\ell\,(H(x)+s_S)$, so $S\circ A_t$ is represented by
$\red(S)+T_t\,(H+s_S)$ with $T_t=\sum_it^{2^i}X^{2^i}$. The products $X^{2^i}X^{2^k}$ ($k\ge1$) have exponent $<q$
except $i=k=n-1$, where $X^q$ reduces to $X$; hence
$\red(S\circ A_t)=\red(S)+T_t(H+s_S)+t^{q/2}u\,(X+X^q)$ and, differentiating,
\[\red(S\circ A_t)'=t\,(H+s_S)+t^{q/2}u .\]
Its linear part $tH$ has kernel $\{0,1\}$, so rank $n-1$, and the constant lies outside
$\im(tH)=t\,\im H$ iff $\Tr(\pi(s_S+t^{q/2-1}u))=1$, i.e.\ iff $\Tr(\pi u\,t^{-1/2})=0$, i.e.\ (squaring inside the
trace) iff $\Tr(\gamma/t)=0$ with $\gamma=\pi^2u^2$. If $\gamma=0$ any $t\ne0$ with $\Tr t=0$ works; otherwise
Lemma~\ref{lem:trace4} with $(\beta,\varepsilon)=(0,0)$ gives one.
\end{proof}

\begin{proposition}[kernel line]\label{prop:outKL}
Let $n\ge5$ be odd, $Q$ a quadratic APN permutation, $A$ invertible, $g=Q\circ A$, $c(\ell)=p\cdot\ell$. Suppose
$w\in W_A=\ker_FN$ has $c(w)\ne0$ and $\rk L_w\ge n-1$. Then for suitable $a\in F^\times$ and $f\in F^\times$ the map
$L=aL_w+\Tr(f\,\cdot)$ is bijective and $L\circ Q\circ A$ is optimal.
\end{proposition}
\begin{proof}
Let $v=m(f)$, so $L_v=\Tr(f\,\cdot)$ and $L=L_{aw+v}$. Since $w\in\ker N$, $K_{aw+v}=aK_w+K_v=K_v=K_f$, which has rank
$n-1$, kernel $\{0,r\}$ and image $\{x:\Tr(rx)=0\}$ (Lemma~\ref{lem:comp}); and $c(aw+v)=a\,c(w)+c_f$. So the constant
is outside the image iff $\Tr(\rho a)=\beta$ with $\rho=r\,c(w)\ne0$ and $\beta=1+\Tr(rc_f)$.

\emph{Bijectivity.} If $L_w$ is bijective, choose any $f\ne0$ and write $\Tr(f\,\cdot)=\Tr(f_1L_w(\cdot))$ with
$f_1\ne0$; then $L=(a\,\mathrm{id}+\Tr(f_1\,\cdot))\circ L_w$, which is bijective iff $a\ne0$ and $\Tr(f_1/a)=0$ (as in
Proposition~\ref{prop:outC}). If $\ker L_w=\{0,z_0\}$ and $\im L_w=\{u:\Tr(\tau u)=0\}$ ($\tau\ne0$), choose $f$ with
$\Tr(fz_0)=1$. A root $y\ne0$ of $L$ with $\Tr(fy)=0$ would have $L_w(y)=0$, so $y=z_0$, contradicting $\Tr(fz_0)=1$.
A root with $\Tr(fy)=1$ needs $L_w(y)=1/a$, which is solvable iff $\Tr(\tau/a)=0$, and then exactly one of the two
preimages $y,y+z_0$ has $\Tr(fy)=1$. So $L$ is bijective iff $a\ne0$ and $\Tr(\tau/a)=1$.

\emph{Existence.} We need $a\ne0$ with $\Tr(\rho a)=\beta$ and $\Tr(\kappa/a)=\varepsilon$, where
$(\kappa,\varepsilon)=(f_1,0)$ or $(\tau,1)$. With $z=\rho a$ and $\gamma=\rho\kappa\ne0$ this is
Lemma~\ref{lem:trace4}. Then $L\circ Q\circ A$ is a quadratic APN permutation whose reduced derivative has rank-$(n-1)$
linear part and constant outside its image.
\end{proof}

\subsection{The two-invariant criterion}\label{ss:twoinv}
Throughout \S\ref{ss:twoinv}--\ref{ss:relabel}, $n\ge5$ is odd, $Q$ is a quadratic APN permutation, $A$ is
invertible, $g=Q\circ A$, $c(\ell)=p\cdot\ell$, and $W,d,U,s,\eta$ are those of $g$ (Section~\ref{sec:prelim}).

\begin{proposition}[the constant invariant]\label{prop:L1eta}
$\eta=1$ iff $p\notin\row_FN$ iff $\rank\binom{p}{N}=\rank N+1$; otherwise $\rank\binom pN=\rank N$.
\end{proposition}
\begin{proof}
For the nondegenerate $F$-bilinear form $u\cdot v$ on $F^n$, $(\ker N)^\perp=\row N$ (the inclusion $\supseteq$ is
immediate and both sides have dimension $\rank N$). So $c|_W\equiv0$ iff $p\in\row N$. The rank of $\binom pN$
exceeds that of $N$ by one exactly when the added row $p$ lies outside $\row N$.
\end{proof}

\begin{lemma}[all rank strata]\label{lem:L1local}
There is $w\in W$ with $c(w)\ne0$ and $\rk L_w\ge n-1$ if and only if $s\le1$ and $\eta=1$. When this holds, for
every $\delta\ne0$ the witness can be chosen with $c(w)=\delta$ and $\ker L_w=U$. No assumption on $\rank N$ is made.
\end{lemma}
\begin{proof}
($\Leftarrow$) $\eta=1$ makes $W\ne0$ and $c|_W$ a nonzero $F$-linear functional; Theorem~\ref{thm:L1lab} gives
$w\in W$ with $c(w)=\delta$ and $\ker L_w=U$, so $\rk L_w=n-s\ge n-1$. ($\Rightarrow$) $c(w)\ne0$ gives $\eta=1$;
every $L_w$ with $w\in W$ kills $U$, so $n-1\le\rk L_w\le n-s$.
\end{proof}

\begin{theorem}[two-invariant criterion]\label{thm:twoinv}
Let $n\ge5$ be odd, $Q$ a quadratic APN permutation and $A$ invertible.
\begin{enumerate}
\item[(i)] There is $w\in W_A$ with $p_A\cdot w\ne0$ and $\rk w\ge n-1$ if and only if $s(A)\le1$ and
$\eta(A)=1$.
\item[(ii)] At the fixed input $A$: CD (some $w\in W_A$ with $L_w$ bijective) holds iff $s(A)=0$;
nonzero-constant CD (some such $w$ with $p_A\cdot w\ne0$) holds iff $s(A)=0$ and $\eta(A)=1$; the kernel-line
hypothesis of Proposition~\ref{prop:outKL} holds iff $s(A)\le1$ and $\eta(A)=1$. CD with $c|_{W_A}\equiv0$ ($s=0$, $\eta=0$) is not
contained in the kernel-line case (example: $Q=x^3$, $A=x^2$, $n=5$). Nothing is asserted here about the existence
of a good $A$.
\item[(iii)] If $s(A)=0$, or $s(A)=1$ and $\eta(A)=1$, then there are invertible $\F_2$-linear maps $L,A'$ such that
$L\circ Q\circ A'$ is optimal, obtained explicitly: via Proposition~\ref{prop:outC} if $(s,\eta)=(0,1)$, via Proposition~\ref{prop:outNF}
(input shear) if $(s,\eta)=(0,0)$, via Proposition~\ref{prop:outKL} if $(s,\eta)=(1,1)$.
\end{enumerate}
\end{theorem}
\begin{proof}
(i) is Lemma~\ref{lem:L1local}. (ii) The first equivalence is Corollary~\ref{cor:L1bij}(b). For the second,
($\Leftarrow$) is Theorem~\ref{thm:L1lab} with $s=0$ and any $\delta\ne0$, and ($\Rightarrow$) holds because a
bijective $L_w$ forces $U=0$ and $c(w)\ne0$ forces $\eta=1$. The third is (i) together with
Corollary~\ref{cor:L1norm}. For the example: $x^3\circ x^2=x^6$ has $d=3$, $U=0$, $s=0$, $\eta=0$
(direct computation; Supplement~S3).

(iii) In branches (1) and (2) below, let $h=L_w\circ Q\circ A$ with $w\in W_A$ and $L_w$ bijective. Then $h$ is linearly
equivalent to $Q$, hence a quadratic APN permutation, and $\red(h)'=c(w)$ is constant because $w\in\ker N$
(\S\ref{ss:data}). In branch (3) $L_w$ has rank $n-1$ and $L_w\circ Q\circ A$ is not a permutation; it is not used there.
\begin{enumerate}
\item $s=0$, $\eta=1$: choose $w\in W_A$ with $\ker L_w=0$ and $c(w)=\delta\ne0$ (Theorem~\ref{thm:L1lab}). Then
$\red(h)'=\delta$, and Proposition~\ref{prop:outC} gives an output update $y\mapsto\sigma y+\Tr(\theta y)$ making $h$
optimal.
\item $s=0$, $\eta=0$: choose $w\in W_A$ with $L_w$ bijective (Theorem~\ref{thm:L1max}); then $\red(h)'=0$ and
Proposition~\ref{prop:outNF} gives $t\ne0$, $\Tr t=0$, with $h\circ A_t$ optimal, $A_t(x)=x+\Tr(tx)$.
\item $s=1$, $\eta=1$: (i) gives $w\in W_A$ with $c(w)\ne0$ and $\rk L_w=n-1$, which is the hypothesis of
Proposition~\ref{prop:outKL}; Proposition~\ref{prop:outKL} gives an optimal $L_{\ell}\circ Q\circ A$ with $\ell=aw+m(f)$,
$L_{m(f)}=\Tr(f\,\cdot)$.
\end{enumerate}
\end{proof}

\begin{remark}[form of the representative]\label{rem:form}
In each branch of (iii) the optimal representative has the form $P=L\circ Q\circ A'$ with $L$ and $A'$ invertible
$\F_2$-linear maps: in (1) $L$ is the output update composed with $L_w$ and $A'=A$; in (2) $A'=A\circ A_t$ with
$A_t(x)=x+\Tr(tx)$, which is bijective because $\Tr t=0$; in (3) $L=L_\ell$ and $A'=A$. In every case $A'$ is bijective
and $P$ is a permutation (optimal representatives are permutations by definition), so the linear map $L$ is bijective
as well. This is also the form of the formal statement (Section~\ref{sec:lean}).
\end{remark}

\begin{remark}
The criterion is sufficient only: it describes when the three constructions above succeed at the given input, not
when an optimal representative $L\circ Q\circ A$ exists for that $A$. Inputs outside the criterion are not needed:
Sections~\ref{sec:layer2}--\ref{sec:small} show that an input satisfying it always exists. Same-input constructions
for some other values of $(s,\eta)$ are collected in Supplement~S2; they are not used.
\end{remark}

\subsection{Output relabeling}\label{ss:relabel}
\begin{proposition}\label{prop:L1relabel}
Let $B$ be an invertible $\F_2$-linear map and $g'=B\circ g$. Then $d,s,\eta$ of $g'$ equal those of $g$, and
$U_{g'}=B(U_g)$.
\end{proposition}
\begin{proof}
$L\mapsto L\circ B$ is a bijection of $\mathrm{Hom}_{\F_2}(F,F)$ that is $F$-linear for the left action. By
\S\ref{ss:data}, $L\in W_{g'}$ (i.e.\ $\red(L\circ B\circ g)'$ is constant) iff $L\circ B\in W_g$; so
$W_{g'}=\{L\circ B^{-1}:L\in W_g\}$, of the same $F$-dimension, with common kernel $B(U_g)$. The constant is
transported unchanged, $c_{g'}(L\circ B^{-1})=c_g(L)$, since both are the constant term of $\red(L\circ g)'$.
\end{proof}

\begin{remark}[scope]\label{rem:L1scope}
Among the linear-equivalence representatives $B\circ Q\circ A$, once the full space $W_A$ is used, re-choosing the
output cannot repair a failing $(s,\eta)$. For a general EA member $B\circ Q\circ A+C$ with $C$ affine, $N$ is
unchanged (it depends only on the quadratic part), hence $W$, $U$, $s$ are unchanged; but $c$ changes by the
$F$-linear functional $\ell\mapsto[X^1]\red(L_\ell\circ C)$, so invariance of $\eta$ is \emph{not} claimed.
Translations are covered by Proposition~\ref{prop:L1Np}. Whether a permutation-preserving linear $C$ can change
$\eta$ is not needed here.
\end{remark}

\section{Global linearization in \texorpdfstring{$\zeta_A$}{zeta\_A}}\label{sec:layer2}

This section proves the linearization used in Section~\ref{sec:mainproof}. Throughout, $g$ is quadratic, i.e.\ every exponent of
$\red(g)$ has binary weight at most $2$; APN and bijectivity are not used in this section.

\subsection{Dickson matrices and the polar tensor}
For $\F_2$-linear $A(x)=\sum_k\alpha_kx^{2^k}$ let $\mathcal D_A[c,j]=\alpha_{j-c}^{2^c}$ (indices mod $n$), so that
$m(Ax)=\mathcal D_A\,m(x)$, and put
\[\zeta_A:=\mathcal D_Ae_0,\qquad(\zeta_A)_c=\alpha_{-c}^{2^c}=[x^1]\,(Ax)^{2^c}.\]
The map $A\mapsto\zeta_A$ is a bijection from $\mathrm{End}_{\F_2}(F)$ onto $F^n$ (it is a coordinatewise
Frobenius twist of $(\alpha_k)_k$). We use three standard facts.
\begin{enumerate}
\item[(F1)] $\mathcal D_A$ is invertible iff $A$ is bijective, and $\det\mathcal D_A\in\F_2$. \emph{Proof.} For an
$\F_2$-basis $x_1,\dots,x_n$ of $F$ the Moore matrix $X=[m(x_1)\cdots m(x_n)]$ is invertible
\cite{Moore}, and $\mathcal D_AX=[m(Ax_1)\cdots m(Ax_n)]$ is invertible iff the $Ax_k$ are $\F_2$-independent.
Squaring every entry of $\mathcal D_A$ permutes rows and columns by the same cyclic shift, so
$\det(\mathcal D_A)^2=\det\mathcal D_A$.
\item[(F2)] For quadratic $\varphi$ with reduced coefficients $c_e(\varphi)$,
$\varphi(x+y)+\varphi(x)+\varphi(y)+\varphi(0)=\sum_{a\ne c}c_{2^a+2^c}(\varphi)\,x^{2^a}y^{2^c}$ (ordered pairs);
constant and weight-one monomials contribute nothing.
\item[(F3)] For quadratic $h$ and $\varphi=L_\ell\circ h$: $c_{1+2^j}(\varphi)=(N_h\ell)_j$ and
$c_1(\varphi)=p_h\cdot\ell$ (definition of $N_h,p_h$, \S\ref{ss:data}).
\end{enumerate}

Let $c_e$ be the reduced coefficients of $g$. Put $\Gamma_{a,c,i}=c_{2^{a-i}+2^{c-i}}^{\,2^i}$ for $a\ne c$,
$\Gamma_{a,a,i}=0$, and for $\xi,\chi\in F^n$ let $\Gamma(\xi,\chi)_i=\sum_{a,c}\xi_a\chi_c\Gamma_{a,c,i}$.

\begin{lemma}\label{lem:gamma}
$\Gamma$ is $F$-bilinear, symmetric and alternating, and for $x,y\in F$, $\ell\in F^n$:
$\Gamma(m(x),m(y))\cdot\ell=L_\ell(B_g(x,y))$, $B_g(x,y)=g(x+y)+g(x)+g(y)+g(0)$.
\end{lemma}
\begin{proof}
Symmetry is $\Gamma_{a,c,i}=\Gamma_{c,a,i}$; then $\Gamma(\xi,\xi)_i=2\sum_{a<c}\xi_a\xi_c\Gamma_{a,c,i}=0$.
By (F2), $B_g(x,y)^{2^i}=\sum_{a\ne c}c_{2^a+2^c}^{2^i}x^{2^{a+i}}y^{2^{c+i}}=\sum_{a\ne c}\Gamma_{a,c,i}x^{2^a}y^{2^c}$
after the shift $(a,c)\mapsto(a-i,c-i)$ and $x^{2^n}=x$. Multiply by $\ell_i$ and sum over $i$.
\end{proof}

Let $G(\zeta)$ be the $n\times n$ matrix with rows $G(\zeta)_c=\Gamma(\zeta,e_c)$, and $\bG_\ell$ the alternating
matrix $\bG_\ell[a,c]=\sum_i\Gamma_{a,c,i}\ell_i$. Then $G$ is linear in $\zeta$, $\bG$ is linear in $\ell$, and
\begin{equation}\label{eq:GGbar}
\chi^\top G(\zeta)\ell=\Gamma(\zeta,\chi)\cdot\ell=\zeta^\top\bG_\ell\chi,\qquad G(\zeta)\ell=\bG_\ell\zeta,\qquad
\zeta^\top G(\zeta)=0 .
\end{equation}

\subsection{The linearization theorem}
\begin{theorem}[global linearization]\label{thm:lin}
Let $A$ be $\F_2$-linear, $h=g\circ A$ and $\zeta=\zeta_A$.
\begin{enumerate}
\item[(a)] $\binom{0}{N_h}=\mathcal D_A^\top G(\zeta)$.
\item[(b)] If $A$ is invertible, $\ker_FN_h=\ker G(\zeta)=\{\ell:\bG_\ell\zeta=0\}$. So $W_A$, $d(A)$, $U_A$ and
$s(A)$ depend on $A$ only through $[\zeta_A]\in\mathbb P^{n-1}(F)$.
\item[(c)] There is $\omega(Z)\in F[Z_0,\dots,Z_{n-1}]^n$, each coordinate zero or homogeneous of degree $n-2$, with
$\mathrm{minor}_{i,k}G(Z)=Z_i\,\omega_k(Z)$ for all $i,k$ ($\mathrm{minor}_{i,k}$ deletes row $i$ and column $k$).
For every $\F_2$-linear $A$ the vector $C(N_h)$ of maximal minors of $N_h$ (column $k$ deleted in coordinate $k$)
equals $\det(\mathcal D_A)\,\omega(\zeta_A)$; for invertible $A$ it equals $\omega(\zeta_A)$.
\item[(d)] (label formula) $p_h=\lambda(\zeta)+\kappa(\zeta)^{2^{n-1}}$ coordinatewise, where
$\lambda_i(Z)=\sum_ac_{2^{a-i}}^{\,2^i}Z_a$ and $\kappa_i(Z)=\sum_{a<c}\Gamma_{a,c,i}^{\,2}Z_{a+1}Z_{c+1}$.
This holds for every $\F_2$-linear $A$, invertible or not.
\item[(e)] $\Phi(Z):=(\lambda(Z)\cdot\omega(Z))^2+\sum_i\kappa_i(Z)\,\omega_i(Z)^2$ is zero or homogeneous of degree
$2n-2$, and for every $\F_2$-linear $A$: $\Phi(\zeta_A)=\bigl(p_{g\circ A}\cdot\omega(\zeta_A)\bigr)^2$.
\end{enumerate}
\end{theorem}
\begin{proof}
(a) $h$ is quadratic and $B_h(w,x)=B_g(Aw,Ax)$. Fix $\ell$, put $\varphi=L_\ell\circ h$, so $B_\varphi=L_\ell\circ B_h$.
By (F2)--(F3) the coefficient of $w^1x^{2^k}$ in $B_\varphi$ (degrees $<q$ in each variable) is $(N_h\ell)_k$ for
$k\ne0$ and $0$ for $k=0$. By Lemma~\ref{lem:gamma},
$B_\varphi(w,x)=\Gamma(\mathcal D_Am(w),\mathcal D_Am(x))\cdot\ell=\sum_{j,k}w^{2^j}x^{2^k}\,\Gamma(\mathcal D_Ae_j,\mathcal D_Ae_k)\cdot\ell$,
whose $w^1x^{2^k}$-coefficient is $\Gamma(\zeta,\mathcal D_Ae_k)\cdot\ell=\sum_c\mathcal D_A[c,k]\,G(\zeta)_c\cdot\ell$,
the $k$-th row of $\mathcal D_A^\top G(\zeta)$ applied to $\ell$.

(b) By (F1) $\mathcal D_A^\top$ is invertible, so $\ker\binom0{N_h}=\ker G(\zeta)$; the second description is
\eqref{eq:GGbar}. $G(t\zeta)=tG(\zeta)$ gives the projective statement.

(c) Over $K=F(Z)$ we have $Z^\top G(Z)=0$ with $Z\ne0$, so $\det G(Z)=0$ and $\mathrm{adj}(G)G=0$: every row of
the adjugate lies in the left kernel of $G(Z)$. If $\rank_KG(Z)=n-1$ this kernel is $KZ^\top$, so row $k$ of
$\mathrm{adj}(G)$ is $\omega_kZ^\top$ with $\omega_k\in K$, i.e.\ $\mathrm{minor}_{i,k}=Z_i\omega_k$ (no signs in
characteristic $2$); if $\rank_KG(Z)<n-1$ all minors vanish and $\omega:=0$. From
$Z_1\,\mathrm{minor}_{0,k}=Z_0\,\mathrm{minor}_{1,k}$ and unique factorization, $Z_0\mid\mathrm{minor}_{0,k}$, so
$\omega_k=\mathrm{minor}_{0,k}/Z_0$ is a polynomial, zero or homogeneous of degree $(n-1)-1$. Let $Y$ be rows
$1,\dots,n-1$ of $\mathcal D_A^\top$, so $N_h=YG(\zeta)$ by (a). Cauchy--Binet gives
$C_k(YG)=\sum_iC_i(Y)\,\mathrm{minor}_{i,k}G=\bigl(C(Y)\cdot\zeta\bigr)\omega_k(\zeta)$, and
$C(Y)\cdot\zeta=\det\binom{\zeta^\top}{Y}=\det\mathcal D_A^\top$ since row $0$ of $\mathcal D_A^\top$ is $\zeta^\top$.
For invertible $A$ this determinant is $1$ by (F1).

(d) $p_{h,i}=[x^1]\red(g(Ax)^{2^i})$. The reduced coefficients of $g^{2^i}$ are $c_{e2^{-i}}^{2^i}$, so
$g(y)^{2^i}=\mathrm{const}+\sum_ac_{2^{a-i}}^{2^i}y^{2^a}+\sum_{a<c}\Gamma_{a,c,i}y^{2^a}y^{2^c}$. Put $y=Ax$. The
$x^1$-coefficient of $(Ax)^{2^a}=\sum_j\mathcal D_A[a,j]x^{2^j}$ is $\zeta_a$. In
$(Ax)^{2^a}(Ax)^{2^c}=\sum_{j,k}\mathcal D_A[a,j]\mathcal D_A[c,k]x^{2^j+2^k}$ the only term reducing to $x^1$ is
$j=k=n-1$: for $j\ne k$ the exponent has weight two, and for $j=k<n-1$ it is $x^{2^{j+1}}\ne x$. Finally
$\mathcal D_A[a,n-1]^2=\alpha_{n-1-a}^{2^{a+1}}=\zeta_{a+1}$, so $\mathcal D_A[a,n-1]=\zeta_{a+1}^{2^{n-1}}$, and
$p_{h,i}=\lambda_i(\zeta)+\sum_{a<c}\Gamma_{a,c,i}(\zeta_{a+1}\zeta_{c+1})^{2^{n-1}}=\lambda_i(\zeta)+\kappa_i(\zeta)^{2^{n-1}}$.
Nothing here uses invertibility.

(e) Homogeneity is clear. For $\zeta\in F^n$, squaring is additive and $\kappa_i(\zeta)^{2^n}=\kappa_i(\zeta)$, so
$\bigl((\lambda(\zeta)+\kappa(\zeta)^{2^{n-1}})\cdot\omega(\zeta)\bigr)^2=(\lambda(\zeta)\cdot\omega(\zeta))^2+\sum_i\kappa_i(\zeta)\omega_i(\zeta)^2=\Phi(\zeta)$;
now use (d).
\end{proof}

\begin{remark}
For invertible $A$, (c) and (e) give $\Phi(\zeta_A)=(\det\binom{p_h}{N_h})^2$, which is nonzero iff $d(A)=1$ and
$\eta(A)=1$ (the reverse implication is not needed below). We prove $\Phi\not\equiv0$ in
Section~\ref{sec:cofactor} by evaluating at Moore points.
\end{remark}

\section{The cofactor dictionary and two nonzero polynomials}\label{sec:cofactor}

From now on $Q$ is a quadratic APN permutation, $n\ge5$ odd, and $G,\bG,\omega,\Phi$ are formed from $Q$.
We separate two ranks: $\rank$ is the rank of a matrix over $F$, and $\rk\ell=\dim_{\F_2}L_\ell(F)$ is the
$\F_2$-rank of the map $L_\ell$.

\begin{lemma}[rank dictionary]\label{lem:rank}
(a) For $\F_2$-linear $A$, $L_{\zeta_A}$ is the adjoint of $A$ for $\langle x,y\rangle=\Tr(xy)$; hence
$\rk\zeta_A=\rank_{\F_2}A$. (b) $\zeta_{A\circ\mu_t}=t\,\zeta_A$ for $\mu_t(x)=tx$, $t\in F^\times$, and
$\rk(t\ell)=\rk\ell$. Hence for $\zeta\ne0$ the set $\{A:\zeta_A\in F^\times\zeta\}$ has exactly $q-1$ elements, all
of rank $\rk\zeta$.
\end{lemma}
\begin{proof}
(a) With $A(x)=\sum_k\alpha_kx^{2^k}$, $\Tr(A(x)y)=\sum_k\Tr\bigl(x(\alpha_ky)^{2^{-k}}\bigr)
=\Tr\bigl(x\sum_c(\alpha_{-c}y)^{2^c}\bigr)=\Tr(x\,L_{\zeta_A}(y))$. An adjoint for a nondegenerate form has the
same rank. (b) $A\circ\mu_t$ has coefficients $\alpha_kt^{2^k}$, so $(\zeta_{A\mu_t})_c=(\alpha_{-c}t^{2^{-c}})^{2^c}=t(\zeta_A)_c$;
$L_{t\ell}=tL_\ell$. The last sentence follows from bijectivity of $A\mapsto\zeta_A$.
\end{proof}

\begin{lemma}[cofactor dictionary]\label{lem:s}
Let $\zeta\in F^n\setminus0$. Then $G(\zeta)\,\omega(\zeta)=0$; $\omega(\zeta)\ne0$ iff $\rank G(\zeta)=n-1$; and then
$\ker G(\zeta)=F\omega(\zeta)$. If $A$ is invertible and $\omega(\zeta_A)\ne0$, then
\[d(A)=1,\qquad W_A=F\,\omega(\zeta_A),\qquad U_A=\ker L_{\omega(\zeta_A)},\qquad s(A)=n-\rk\omega(\zeta_A).\]
\end{lemma}
\begin{proof}
By \eqref{eq:GGbar} $\zeta^\top G(\zeta)=0$, so $\det G(\zeta)=0$ and $G\,\mathrm{adj}(G)=0$. Column $i$ of
$\mathrm{adj}(G)$ is $(\mathrm{minor}_{i,k}G)_k=\zeta_i\,\omega(\zeta)$ (Theorem~\ref{thm:lin}(c)); choosing
$\zeta_i\ne0$ gives $G\omega=0$. The $(n-1)$-minors are the $\zeta_i\omega_k(\zeta)$, and not all $\zeta_i$ vanish, so
some minor is nonzero iff $\omega(\zeta)\ne0$; and $\rank G\le n-1$. For invertible $A$, Theorem~\ref{thm:lin}(b)
gives $W_A=\ker G(\zeta_A)=F\omega(\zeta_A)$, so $d(A)=1$. Since $L_{t\ell}=tL_\ell$, the common kernel is
$\ker L_{\omega(\zeta_A)}$, of dimension $n-\rk\omega(\zeta_A)$.
\end{proof}

\begin{lemma}[$\Phi$ at invertible inputs]\label{lem:phieta}
If $A$ is invertible and $\Phi(\zeta_A)\ne0$, then $d(A)=1$ and $\eta(A)=1$.
\end{lemma}
\begin{proof}
By Theorem~\ref{thm:lin}(e), $\Phi(\zeta_A)=(p_A\cdot\omega(\zeta_A))^2$, so $\omega(\zeta_A)\ne0$ and
$p_A\cdot\omega(\zeta_A)\ne0$. By Lemma~\ref{lem:s}, $d(A)=1$ and $W_A=F\omega(\zeta_A)$, on which
$\ell\mapsto p_A\cdot\ell$ is not identically zero; that is $\eta(A)=1$.
\end{proof}

\subsection{Moore points}
For $x\in F^\times$ put $m(x)=(x^{2^i})_i$. Then $m(x)=\zeta_{A_x}$ for the rank-one map $A_x(y)=x\Tr(y)$
(all $\alpha_k=x$). Moore points are used only to show that two polynomials are nonzero; they are not inputs of
the construction, and no good input is assumed.

We use two facts, labeled (AB) (almost bent) and (MB) (Moore basis): (AB) for $a\in F^\times$ the $\F_2$-bilinear form $b_a(x,z)=\Tr(aB_Q(x,z))$ has a radical of dimension
exactly $1$ (the almost-bent property of quadratic APN functions in odd dimension~\cite{CCZ98}, proved in
Lemma~\ref{lem:L1apn}); (MB) if $x_1,\dots,x_n$ is an $\F_2$-basis of $F$, the Moore vectors
$m(x_1),\dots,m(x_n)$ form an $F$-basis of $F^n$ \cite{Moore}.

\begin{lemma}[kernel at a Moore point]\label{lem:moorekernel}
Let $x\in F^\times$. The image $H_x=\{B_Q(x,z):z\in F\}$ is an $\F_2$-hyperplane $\{y:\Tr(a_xy)=0\}$ with
$a_x\ne0$, $\ker G(m(x))=F\,m(a_x)$, and $\omega(m(x))=t_x\,m(a_x)$ with $t_x\in F^\times$.
\end{lemma}
\begin{proof}
$z\mapsto B_Q(x,z)$ is $\F_2$-linear and, by APN, has kernel $\{0,x\}$; so $H_x$ is a hyperplane. By
\eqref{eq:GGbar}, $\ell\in\ker G(m(x))$ iff $\Gamma(m(x),\chi)\cdot\ell=0$ for all $\chi\in F^n$; by (MB) it
suffices to take $\chi=m(z)$, and by Lemma~\ref{lem:gamma} the condition reads $L_\ell(B_Q(x,z))=0$ for all $z$,
i.e.\ $L_\ell$ vanishes on $H_x$. The $\F_2$-linear maps vanishing on $H_x$ are $y\mapsto b\Tr(a_xy)=L_{b\,m(a_x)}(y)$,
$b\in F$. So $\ker G(m(x))=F\,m(a_x)$ is one-dimensional, $\rank G(m(x))=n-1$, and Lemma~\ref{lem:s} gives
$\omega(m(x))\in\ker G(m(x))\setminus0$.
\end{proof}

\begin{proposition}[$\Phi\not\equiv0$]\label{prop:phinonzero}
For every $x\in F^\times$, $\Phi(m(x))=t_x^2\ne0$. In particular $\Phi\ne0$, and $\Phi$ has total degree at most $2n-2$.
\end{proposition}
\begin{proof}
Let $h=Q\circ A_x$. Since $\Tr(y)\in\F_2$, $h(y)=Q(0)+e\,\Tr(y)$ with $e=Q(x)+Q(0)$, and
$h^{2^i}=Q(0)^{2^i}+e^{2^i}\Tr$; the $X^1$-coefficient of $\red(\Tr)$ is $1$, so $p_h=m(e)$. By
Theorem~\ref{thm:lin}(e) and Lemma~\ref{lem:moorekernel},
$\Phi(m(x))=(m(e)\cdot t_xm(a_x))^2=t_x^2\Tr(a_xe)^2$. If $e\in H_x$, then $Q(x)+Q(0)=Q(x+z)+Q(x)+Q(z)+Q(0)$ for
some $z$, i.e.\ $Q(x+z)=Q(z)$, and injectivity of $Q$ gives $x=0$, a contradiction. So $\Tr(a_xe)=1$.
\end{proof}

\begin{proposition}[the adjugate polynomial $R$]\label{prop:R}
There is a nonzero homogeneous $R\in F[Z]$ of degree $(n-1)(n-2)$ such that for
every $\zeta\in F^n\setminus0$ with $R(\zeta)\ne0$: $\omega(\zeta)\ne0$, $\rank\bG_{\omega(\zeta)}=n-1$ and
$\ker\bG_{\omega(\zeta)}=F\zeta$.
\end{proposition}
\begin{proof}
Fix $x\in F^\times$ and $a=a_x$. For $\xi,\chi\in F^n$, $\xi^\top\bG_{m(a)}\chi=\Gamma(\xi,\chi)\cdot m(a)$, which on
Moore vectors equals $L_{m(a)}(B_Q(x',z))=\Tr(aB_Q(x',z))=b_a(x',z)$. In the $F$-basis $m(x_1),\dots,m(x_n)$ of (MB)
the Gram matrix of $\bG_{m(a)}$ is the $\F_2$-matrix $(b_a(x_j,x_k))_{j,k}$, whose rank is $n-1$ by (AB); rank does
not change under congruence by an invertible matrix or under extension of scalars. So $\rank\bG_{m(a)}=n-1$, and
by Lemma~\ref{lem:moorekernel} and linearity of $\ell\mapsto\bG_\ell$,
$\bG_{\omega(m(x))}=t_x\bG_{m(a)}$ has rank $n-1$. Its adjugate is then nonzero: choose $(i_0,k_0)$ with
$\mathrm{minor}_{i_0,k_0}\bG_{\omega(m(x))}\ne0$ and put $R(Z)=\mathrm{minor}_{i_0,k_0}\bigl(\bG_{\omega(Z)}\bigr)$.
The entries of $\bG_{\omega(Z)}$ are $F$-linear combinations of the $\omega_i(Z)$, hence zero or homogeneous of
degree $n-2$, so $R$ is zero or homogeneous of degree $(n-1)(n-2)$; since $R(m(x))\ne0$, $R$ is nonzero and hence
homogeneous of that degree.

Let $\zeta\ne0$ with $R(\zeta)\ne0$. Then $\bG_{\omega(\zeta)}\ne0$, so $\omega(\zeta)\ne0$, and
$\rank\bG_{\omega(\zeta)}\ge n-1$. By \eqref{eq:GGbar} and Lemma~\ref{lem:s},
$\bG_{\omega(\zeta)}\zeta=G(\zeta)\omega(\zeta)=0$ with $\zeta\ne0$, so the rank is exactly $n-1$ and the kernel is
$F\zeta$.
\end{proof}

\begin{lemma}[projective injectivity]\label{lem:inj}
The map $[\zeta]\mapsto[\omega(\zeta)]$ is injective on $\{[\zeta]\in\mathbb P^{n-1}(F):R(\zeta)\ne0\}$.
\end{lemma}
\begin{proof}
If $R(\zeta)R(\zeta')\ne0$ and $\omega(\zeta')=t\,\omega(\zeta)$ with $t\ne0$, then $\bG_{\omega(\zeta')}=t\bG_{\omega(\zeta)}$,
so $\zeta'\in\ker\bG_{\omega(\zeta')}=\ker\bG_{\omega(\zeta)}=F\zeta$ by Proposition~\ref{prop:R}.
\end{proof}

\noindent The map $\zeta\mapsto\omega(\zeta)$ itself is not injective on vectors: $\omega(t\zeta)=t^{n-2}\omega(\zeta)$,
and $\gcd(n-2,q-1)$ can exceed $1$ (e.g.\ $n=9$: $7\mid511$). The count below therefore works with projective
classes only.

\section{Universal input existence for \texorpdfstring{$n\ge9$}{n >= 9} and proof of the Main Theorem}\label{sec:mainproof}

The exceptional set of inputs is cut out by the polynomial $\Phi R$, of degree at most
\[E_n:=(2n-2)+(n-1)(n-2)=n(n-1).\]
A Pfaffian variant gives the smaller degree $D_n=(n^2+n-4)/2$ and a larger margin; it is recorded in
Appendix~\ref{app:pfaff} and is not needed.

Let $\rho_n(r)=\#\{M\in M_n(\F_2):\rank M=r\}$ and $\rho_n(\le r)=\sum_{j\le r}\rho_n(j)$; $\rho_n(\le n-2)$ counts the
zero matrix. Let
\[\mathrm{Fail}(Q)=\{A\in\GL(F):\ \text{neither }s(A)=0\text{ nor }(s(A)=1\text{ and }\eta(A)=1)\}.\]

\begin{theorem}[counting criterion]\label{thm:count}
For every quadratic APN permutation $Q$ of $F$, $n$ odd, $n\ge5$,
\[|\mathrm{Fail}(Q)|\ \le\ \bigl(\rho_n(\le n-2)-1\bigr)+\#\{\zeta\in F^n:\Phi(\zeta)R(\zeta)=0\}
\ \le\ \rho_n(\le n-2)-1+E_n\,q^{n-1}.\]
\end{theorem}
\begin{proof}
Split $\mathrm{Fail}=\mathrm{Fail}_1\sqcup\mathrm{Fail}_2$ according to whether $\Phi R(\zeta_A)=0$. Since
$A\mapsto\zeta_A$ is injective, $|\mathrm{Fail}_1|\le\#\{\zeta:\Phi R(\zeta)=0\}$.

Let $A\in\mathrm{Fail}_2$ and $\zeta=\zeta_A$ ($\ne0$). By Lemma~\ref{lem:phieta}, $d(A)=1$ and $\eta(A)=1$; by
Lemma~\ref{lem:s}, $s(A)=n-\rk\omega(\zeta)$. So $A$ fails iff $\rk\omega(\zeta)\le n-2$. Membership in
$\mathrm{Fail}_2$ depends only on $[\zeta_A]$: $\Phi R$ is homogeneous, $\omega(t\zeta)=t^{n-2}\omega(\zeta)$ has the
same $\rk$, and by Lemma~\ref{lem:rank}(b) the $q-1$ maps $A'$ with $\zeta_{A'}\in F^\times\zeta$ are invertible
(they are $A\circ\mu_t$). Hence $|\mathrm{Fail}_2|=(q-1)\,\#\mathcal S$ with
$\mathcal S=\{[\zeta_A]:A\in\mathrm{Fail}_2\}$. By Lemma~\ref{lem:inj}, $[\zeta]\mapsto[\omega(\zeta)]$ embeds
$\mathcal S$ into $\{[w]\in\mathbb P^{n-1}(F):\rk w\le n-2\}$. Since $\ell\mapsto L_\ell$ is a bijection
$F^n\to\mathrm{End}_{\F_2}(F)\cong M_n(\F_2)$ and $\rk$ is constant on $F^\times$-classes, that set has
$(\rho_n(\le n-2)-1)/(q-1)$ elements (the $-1$ removes the zero matrix). So $|\mathrm{Fail}_2|\le\rho_n(\le n-2)-1$.

Finally $\Phi R\ne0$ (Propositions~\ref{prop:phinonzero}, \ref{prop:R}; $F[Z]$ is a domain) and
$\deg\Phi R\le(2n-2)+(n-1)(n-2)=E_n$. A nonzero polynomial of total degree $D$ in $m$ variables over a field $K$ has
at most $D|S|^{m-1}$ zeros in $S^m$ for every finite $S\subseteq K$ (Ore, Schwartz--Zippel
\cite{Ore,DeMilloLipton,Zippel,Schwartz}); with $S=F$, $m=n$ this is $E_nq^{n-1}$.
\end{proof}

\begin{lemma}[the rank strata]\label{lem:rankcount}
Put $\varpi_n=\prod_{k=1}^n(1-2^{-k})$ and $\Lambda_n=4\varpi_n(1-2^{-n-1})-1$. Then $\rho_n(n)=|\GL_n(\F_2)|=\varpi_nq^n$,
$\rho_n(n-1)=2(1-2^{-n})\,\varpi_nq^n$, and $\rho_n(n)-\rho_n(\le n-2)=\Lambda_nq^n$.
\end{lemma}
\begin{proof}
$|\GL_n(\F_2)|=\prod_{i<n}(2^n-2^i)=\varpi_n2^{n^2}$ and $2^{n^2}=q^n$. The number of matrices of rank $r$ is
$\prod_{i<r}(2^n-2^i)^2/\prod_{i<r}(2^r-2^i)$. For $r=n-1$: $\prod_{i<n-1}(2^n-2^i)=|\GL_n|/2^{n-1}$ and
$\prod_{i<n-1}(2^{n-1}-2^i)=|\GL_{n-1}|$, with $|\GL_n|/|\GL_{n-1}|=2^{n-1}(2^n-1)$; so
$\rho_n(n-1)=|\GL_n|\,(2^n-1)/2^{n-1}=2(1-2^{-n})|\GL_n|$. Since $\sum_r\rho_n(r)=q^n$,
$\rho_n(n)-\rho_n(\le n-2)=2\rho_n(n)+\rho_n(n-1)-q^n=\bigl(\varpi_n(4-2^{1-n})-1\bigr)q^n=\Lambda_nq^n$.
\end{proof}

\begin{theorem}[universal input existence]\label{thm:universal}
Let $n$ be odd, $n\ge9$, and $Q$ a quadratic APN permutation of $\F_{2^n}$. Then some invertible $\F_2$-linear $A$
has $s(A)=0$, or $s(A)=1$ and $\eta(A)=1$. More precisely,
\[|\GL_n(\F_2)|-|\mathrm{Fail}(Q)|\ \ge\ \bigl(\Lambda_n-n(n-1)2^{-n}\bigr)q^n+1\ >\ 0 .\]
\end{theorem}
\begin{proof}
The displayed inequality is Theorem~\ref{thm:count} with Lemma~\ref{lem:rankcount}. It remains to show
$\Lambda_n>n(n-1)2^{-n}$ for odd $n\ge9$. Since $\prod_{j\ge4}(1-2^{-j})\ge1-\sum_{j\ge4}2^{-j}=\tfrac78$,
\[\varpi_n\ \ge\ \prod_{j\ge1}(1-2^{-j})\ \ge\ \tfrac12\cdot\tfrac34\cdot\tfrac78\cdot\tfrac78=\tfrac{147}{512},\]
so for $n\ge9$
\[\Lambda_n\ \ge\ 4\cdot\tfrac{147}{512}\bigl(1-2^{-10}\bigr)-1=\tfrac{77236}{524288}\ >\ \tfrac{73728}{524288}=\tfrac{72}{512}.\]
The sequence $n(n-1)2^{-n}$ is strictly decreasing for $n\ge4$ (the ratio of consecutive terms is
$(n+1)/(2(n-1))<1$), and equals $72/512$ at $n=9$. Hence $\Lambda_n>n(n-1)2^{-n}$ for every $n\ge9$.
\end{proof}

\begin{remark}[numbers; denominators]\label{rem:numbers}
The margins $\Lambda_n-E_n/q$ (denominator $q^n=2^{n^2}$, the number of all $n\times n$ binary matrices) are
$0.015656$ ($n=9$), $0.101724$ ($n=11$), $0.136180$ ($n=13$), and tend to $0.155152$. Relative to $|\GL_n(\F_2)|$ the
lower bound for the fraction of successful inputs is $0.0541$ at $n=9$ and $0.352$ at $n=11$
(Corollary~\ref{cor:proportion}). On the Pfaffian route
(Appendix~\ref{app:pfaff}) the margins are $0.072297$, $0.124184$, $0.144359$. The values for finite $n$ are exact
rational computations; the limit $0.155152\ldots$ is $4\varpi_\infty-1$ with $\varpi_\infty=\prod_{k\ge1}(1-2^{-k})$,
enclosed in a rational interval by $\varpi_{40}(1-2^{-40})\le\varpi_\infty\le\varpi_{40}$ (Supplement~S3). They are
illustrations, and the proof of
Theorem~\ref{thm:universal} for all $n\ge9$ is the inequality above, not these values. For $n=5,7$ the margins are
negative ($-0.452$ and $-0.168$ on this route) and Section~\ref{sec:small} is used instead. The formal proof
(Section~\ref{sec:lean}) contains the inequality $(\rho_n(\le n-2)-1)+n(n-1)q^{n-1}<|\GL_n(\F_2)|$ for all odd $n\ge9$,
not these numerical values.
\end{remark}

\begin{remark}[what the proof uses]\label{rem:circ}
Theorem~\ref{thm:count} uses Theorem~\ref{thm:lin}, the almost-bent property (AB), the Moore determinant (MB),
injectivity of $Q$ (in Proposition~\ref{prop:phinonzero}), elementary linear algebra and the Schwartz--Zippel bound.
None of these assumes a good input, an optimal output, or a classification. The argument is a pigeonhole, not an
average: the rank-$\le n-2$ stratum of $M_n(\F_2)$ ($\approx0.1336\,q^n$) is smaller than the invertible stratum
($\approx0.2888\,q^n$), and the only loss is the hypersurface $\{\Phi R=0\}$ of degree $n(n-1)$.
\end{remark}

\begin{corollary}[proportion of good inputs]\label{cor:proportion}
Let $n\ge9$ be odd and $Q$ a quadratic APN permutation of $\F_{2^n}$. Among the invertible $\F_2$-linear $A$, the
proportion of those with $\Phi R(\zeta_A)\ne0$ and $s(A)\le1$ (for which $d(A)=1$, $\eta(A)=1$, and the criterion of
Theorem~\ref{thm:twoinv} holds) is at least
\[\pi_n=4\bigl(1-2^{-n-1}\bigr)-\frac{1+n(n-1)2^{-n}}{\varpi_n}.\]
In particular $\pi_9>0.054$, $\pi_{11}>0.352$, $\pi_n>0.35$ for every odd $n\ge11$ (so $\pi_n>0.054$ for every odd
$n\ge9$), $\pi_n>1/2$ for every odd $n\ge15$, and
$\pi_n\to4-1/\varpi_\infty=0.5372\ldots$, where $\varpi_\infty=\prod_{k\ge1}(1-2^{-k})$. These are values of the
lower bound $\pi_n$, not of the proportion itself.
\end{corollary}
\begin{proof}
The proof of Theorem~\ref{thm:count} shows that every invertible $A$ is of this kind unless $\Phi R(\zeta_A)=0$ (at
most $E_nq^{n-1}$ inputs) or $\rk\omega(\zeta_A)\le n-2$ (at most $\rho_n(\le n-2)-1$ inputs); by Lemma~\ref{lem:phieta}
and Lemma~\ref{lem:s} the others have $d(A)=1$, $\eta(A)=1$ and $s(A)=n-\rk\omega(\zeta_A)\le1$. Dividing the bound of
Theorem~\ref{thm:universal} by $|\GL_n(\F_2)|=\varpi_nq^n$ and using $\Lambda_n/\varpi_n=4(1-2^{-n-1})-1/\varpi_n$
gives $\pi_n$. For odd $n\ge15$: $1-2^{-n-1}\ge1-2^{-16}$; $n(n-1)2^{-n}\le210\cdot2^{-15}$, since the sequence
decreases for $n\ge4$; and $\varpi_n\ge\varpi_\infty\ge\varpi_{10}\prod_{k>10}(1-2^{-k})\ge\varpi_{10}(1-2^{-10})$.
Hence $\pi_n\ge4(1-2^{-16})-(1+210\cdot2^{-15})/\bigl(\varpi_{10}(1-2^{-10})\bigr)>0.5149$. Likewise, for odd $n\ge11$,
$1-2^{-n-1}\ge1-2^{-12}$, $n(n-1)2^{-n}\le110\cdot2^{-11}$ and the same bound for $\varpi_n$ give
$\pi_n\ge4(1-2^{-12})-(1+110\cdot2^{-11})/\bigl(\varpi_{10}(1-2^{-10})\bigr)>0.35$. These numbers and the values at
$n=9,11$ are exact rational computations; the limit is enclosed in a rational interval using
$\varpi_{40}(1-2^{-40})\le\varpi_\infty\le\varpi_{40}$ (Supplement~S3).
\end{proof}

\begin{remark}[a randomized construction]\label{rem:algorithm}
Corollary~\ref{cor:proportion} makes the Main Theorem for odd $n\ge9$ effective. Draw inputs $A$ independently and uniformly from $\GL_n(\F_2)$;
compute the reduced coefficients of $Q\circ A$ from those of $Q$ (at most $1+n+\binom n2$ are nonzero) with $O(n^4)$
operations in $F$, then $G(\zeta_A)$, $W_A=\ker G(\zeta_A)$ (Theorem~\ref{thm:lin}), $U_A$, $s(A)$ and $\eta(A)$ by
linear algebra; repeat until $d(A)=1$, $\eta(A)=1$ and $s(A)\le1$. By Corollary~\ref{cor:proportion} each draw succeeds with probability at
least $\pi_n$, so the expected number of draws is at most $1/\pi_n$ (fewer than $2$ for $n\ge15$). Then $W_A=F\omega(\zeta_A)$ has $q$ elements, and the constructions of
Propositions~\ref{prop:outC}--\ref{prop:outKL} need parameters in $F$ with explicit trace conditions, found by
testing at most $q^2$ candidates. Every step is checkable, and the output $L\circ Q\circ A'$ can be verified
directly. The expected number of operations, including the construction and the final check, is polynomial in
$q=2^n$, the size of the look-up table of $Q$; it is not polynomial in $n$.
\end{remark}

\begin{proof}[Proof of the Main Theorem]
Let $Q$ be a quadratic APN permutation of $\F_{2^n}$, $n$ odd, $n\ge5$. If $n\ge9$, Theorem~\ref{thm:universal} gives
an invertible $A$ with $s(A)=0$, or $s(A)=1$ and $\eta(A)=1$. If $n\in\{5,7\}$, Proposition~\ref{prop:small} gives an
invertible $A$ with $s(A)=0$. In both cases Theorem~\ref{thm:twoinv}(iii) gives an optimal $P=L\circ Q\circ A'$
with $L,A'$ invertible $\F_2$-linear (Remark~\ref{rem:form}). For $n\in\{5,7\}$ the argument depends on the
classifications cited in Section~\ref{sec:small}.
\end{proof}

\section{Dimensions three, five and seven}\label{sec:small}

For $n=5,7$ the degree bound of Theorem~\ref{thm:count} does not close (Remark~\ref{rem:numbers}). This is not a
defect of the general argument; these two dimensions are treated with the published classifications and
finite certificates. For $n=3$ the statement is false (Proposition~\ref{prop:nthree}). Nothing in this section is
formalized.

\begin{lemma}[permutation system]\label{lem:permsys}
Let $S$ be a quadratic APN function of $F$ and, for $a\ne0$, let $\nu_a\ne0$ span the orthogonal
complement (for the standard inner product on $\F_2^n$, with $F$ identified with $\F_2^n$ through a fixed basis) of the image of $B_S(a,\cdot)$. For a linear
map $L$, $S+L$ is a permutation iff $\langle\nu_a,S(a)+S(0)+La\rangle=1$ for all $a\ne0$. An EA class of
quadratic APN functions containing $S$ contains a permutation iff this binary linear system in the $n^2$
entries of $L$ is consistent.
\end{lemma}
\begin{proof}
$(S+L)(x+a)+(S+L)(x)=B_S(a,x)+S(a)+S(0)+La$, and the image of $B_S(a,\cdot)$ is a hyperplane; so $S+L$ is
injective iff $S(a)+S(0)+La$ lies outside that hyperplane for every $a\ne0$. If $g=A_1\circ S\circ A_2+L_3$
with $A_1(y)=My+c$ and $A_1,A_2$ affine bijections, then $g$ is a bijection iff
$z\mapsto S(z)+M^{-1}L_3(A_2^{-1}z)$ is, and the latter is $S$ plus an affine map; constants do not affect
bijectivity.
\end{proof}

\begin{lemma}[transport of $s$ along EA-equivalence]\label{lem:transport}
Let $S,Q$ be quadratic and EA-equivalent, $Q=C\circ S\circ B+L+\gamma$ with $B(x)=B_0x+\beta$, $C(y)=C_0y+\gamma'$,
$B_0,C_0$ invertible $\F_2$-linear, $L$ linear and $\beta,\gamma,\gamma'\in F$. If $A$ is invertible, then
$s(Q\circ B_0^{-1}A)=s(S\circ A)$.
\end{lemma}
\begin{proof}
With $A'=B_0^{-1}A$: $(Q\circ A')(x)=C_0\bigl((S\circ A)(x+A^{-1}\beta)\bigr)+(\text{affine in }x)$. The matrix $N$ of a
quadratic map depends only on its reduced coefficients of weight two, so added affine terms do not change $N$, hence
not $W,U,s$. The translation $x\mapsto x+A^{-1}\beta$ does not change $N$ (Proposition~\ref{prop:L1Np}),
and the output map $C_0$ does not change $s$ (Proposition~\ref{prop:L1relabel}).
Invariance of $\eta$ under added linear terms is \emph{not} claimed (Remark~\ref{rem:L1scope}), and is not needed:
we only transport inputs with $s=0$.
\end{proof}

\begin{theorem}\label{thm:57}
For $n\in\{5,7\}$ the Main Theorem holds: every quadratic APN permutation $Q$ of $\F_{2^n}$ has invertible
$\F_2$-linear $L,A$ with $L\circ Q\circ A$ optimal. For $n=7$, among the $488$ EA classes of quadratic APN functions exactly three contain a permutation (classes
$0,1,2$ of the tabulation described below); each of the other $485$ has an explicit XOR inconsistency certificate
of $7$ to $25$ rows.
\end{theorem}

\begin{proposition}[input witnesses]\label{prop:small}
For every quadratic APN permutation $Q$ of $\F_{2^5}$ or $\F_{2^7}$ there is an invertible $\F_2$-linear $A$ with
$s(Q\circ A)=0$.
\end{proposition}

\begin{proof}[Proof of Theorem~\ref{thm:57} and Proposition~\ref{prop:small}]
\emph{$n=5$.} By Brinkmann--Leander~\cite[Theorem~6 and \S5]{BrinkmannLeander} every APN function of $\F_{2^5}$ is CCZ-equivalent to a
power function, and by Yoshiara~\cite[Theorem~2]{Yoshiara16} a quadratic APN function is CCZ-equivalent to a
power APN function iff it is EA-equivalent to a Gold function. So every quadratic APN permutation $Q$ of $\F_{2^5}$
is EA-equivalent to $x^{2^k+1}$ for some $k$ with $\gcd(k,5)=1$; up to the linear equivalence
$x^{2^{n-k}+1}=(x^{2^k+1})^{2^{n-k}}$ we may take $k\in\{1,2\}$. For $f=x^3$ and $f=x^5$ explicit invertible $A$ with
$(d,s,\eta)=(1,0,1)$ for $f\circ A$ are listed in Supplement~S3 and were verified by direct interpolation of
$f\circ A$. Lemma~\ref{lem:transport} gives an
invertible input with $s=0$ for $Q$.

\emph{$n=7$.} The list of Kalgin and Idrisova is complete up to CCZ-equivalence~\cite[\S5]{KI}; we use it as tabulated
in the sboxU library~\cite{sboxU} at commit \texttt{826c4292}; the extracted table and its SHA-256 digest are
part of the supplementary material (Supplement~S3). By Yoshiara~\cite[Theorem~1]{Yoshiara12} CCZ- and EA-equivalence coincide for quadratic APN
functions, so the $488$ tabulated functions represent all EA classes. For each class the system of
Lemma~\ref{lem:permsys} was solved: it is consistent for exactly the classes $0,1,2$ (whose tabulated
representatives are themselves permutations), and for each of the other $485$ an explicit subset of $7$ to $25$
equations sums to $0=1$. Consistency was decided by Gaussian elimination over $\F_2$ and, in an independent computation,
with the \texttt{IsConsistent} routine of Magma~\cite{Magma} (version V2.29-7, September 2026); the $485$ certificates were re-verified by independently written code
that rebuilds each equation from the table (Supplement~S3). For each of the classes $0,1,2$ an invertible $A$ with $(d,s,\eta)=(1,0,1)$ for $S\circ A$, $S$ the tabulated
representative, was found and verified by direct interpolation (Supplement~S3). The table is read as a map on
$\F_{2^7}=\F_2[x]/(x^7+x+1)$ in the polynomial basis; another identification changes the map only by linear
equivalence, along which Lemma~\ref{lem:transport} transports an invertible input with $s=0$ (the input is changed
together with the map). Lemma~\ref{lem:transport} again transports the input
to every $Q$ in these classes.

In both dimensions Theorem~\ref{thm:twoinv}(iii) then gives an optimal representative.
\end{proof}

\begin{remark}[what this rests on]\label{rem:small}
The proof for $n=7$ does not identify the classes $0,1,2$ with Gold classes; the witnesses are computed on the
tabulated representatives themselves. The count ``exactly three of the $488$ classes contain a permutation'' is our
computation and is not taken from~\cite{KI}. The proof rests on: Brinkmann--Leander and
Yoshiara~2016 ($n=5$); Kalgin--Idrisova and Yoshiara~2012 ($n=7$), together with the fidelity of the sboxU
tabulation of the Kalgin--Idrisova list; the $485$ XOR certificates; and the five input witnesses ($x^3,x^5$ at
$n=5$; classes $0,1,2$ at $n=7$). The classifications are cited, not re-proved, and are not formalized.

The dependence on one tabulation is removed by a second computation. The classification~\cite{BLLPR26}, obtained with
a different algorithm, lists $2$ EA classes of quadratic APN functions at $n=5$ and $488$ at $n=7$. With independently
written code (standard library only, no code shared with the computation above; Supplement~S3) we solved the system
of Lemma~\ref{lem:permsys} for each of its representatives: both classes at $n=5$ and exactly $3$ of the $488$ at
$n=7$ contain a permutation. For each of these $5$ classes an input $A$ with $s=0$ was found, an optimal
$L\circ Q\circ A'$ was built by Theorem~\ref{thm:twoinv}(iii), and optimality was checked directly (bijective, APN,
derivative without zero, all nonempty fibres of size two). So Theorem~\ref{thm:57} holds with either tabulation as
input.
\end{remark}

\begin{proposition}[$n=3$]\label{prop:nthree}
No quadratic APN permutation $P$ of $\F_8$ has a reduced representative whose formal derivative is free of zeros on
$\F_8$. Hence the Main Theorem and Corollary~\ref{cor:bs} fail for $n=3$, and the Reduced Critical-Point Conjecture
holds for quadratic APN permutations of $\F_8$.
\end{proposition}
\begin{proof}
The reduced exponents of a quadratic map of $\F_8$ lie in $\{0,1,2,4,3,5,6\}$, so there are $8^7$ quadratic maps, and
$\red(P)'=c_1+c_3X^2+c_5X^4$. Exactly $10752$ of them are APN permutations, and none of these has $\red(P)'$ free of
zeros. This was checked by exhaustive enumeration in two independent implementations, a standard-library Python
program and a Magma program with its own field and polynomial arithmetic; both give the same counts (the Magma
program fixes the constant term, which affects neither bijectivity, nor the APN property, nor $\red(P)'$, and
multiplies by $8$; Supplement~S3).
\end{proof}
The same exhaustive conclusion, for all permutations of $\F_8$, was reported earlier by Stone~\cite{Stone-cx}; our two
enumerations are an independent re-verification, included to make the boundary case explicit and reproducible within
the verification package of this paper. The exclusion of $n=3$ is therefore not an artefact of the method; Lemma~\ref{lem:trace4} also fails for $n=3$.

\section{Consequences}\label{sec:cons}

\begin{lemma}[{classical, cf.~\cite{McDonald}; \cite[Theorem~2]{RS}; \cite[Theorem~2.1]{BS-lift}}]\label{lem:grcrit}
Let $\widehat f\in\GR(2^k,n)[X]$, $k\ge2$, reduce to $f\in F[X]$. Then $\widehat f$ permutes
$\GR(2^k,n)$ if and only if $f$ permutes $F$ and $f'$ has no zero on $F$.
\end{lemma}

\begin{proof}[Proof of Corollary~\ref{cor:bs}]
Let $Q$ be a quadratic APN permutation of $\F_{2^n}$, $n$ odd, $n\ge5$, and let $P=L\circ Q\circ A$ be optimal as
in the Main Theorem. Linear equivalence preserves the APN property, so $P$ is an APN permutation of $\F_{2^n}$, and
the formal derivative of its reduced representative has no zero on $\F_{2^n}$. This contradicts
\cite[Conjecture~2]{BS-lift} for this $n$. By Lemma~\ref{lem:grcrit} this is equivalent to saying that, for every $k\ge2$, every lift
of $\red(P)$ permutes $\GR(2^k,n)$, which is the reduced form of the lifting statement discussed in
\cite[Theorem~2.2]{BS-lift}.
\end{proof}

For the Gold classes, closed-form counterexamples are also given by Theorem~\ref{thm:gold} in
Appendix~\ref{app:gold}, and further families are collected in Supplement~S4; neither is needed for
Corollary~\ref{cor:bs}.

\begin{proof}[Proof of Corollary~\ref{cor:ring}]
The argument uses only that $P$ is optimal. Differential uniformity on the ring is $\Delta^R$ of
Section~\ref{sec:prelim}, defined with the difference $\widehat P(x+a)-\widehat P(x)$ as in~\cite[Section~3]{RS};
on $F$ it is the usual $\delta_P$. Let $\widehat P\in\GR(2^k,n)[X]$, $k\ge2$, be any polynomial whose coefficientwise
reduction is $\red(P)$; it need not be a reduced polynomial, and only its reduction is used.
\emph{Permutation.} $\red(P)$ permutes $F$ and $\red(P)'$ has no zero on $F$; the reduction of
$\widehat P'$ is $\red(P)'$, so Lemma~\ref{lem:grcrit} applies.
\emph{Unit directions.} Let $a$ be a unit, with residue $\bar a\ne0$. If $\widehat P(x+a)-\widehat P(x)=b$, then
reducing modulo $2$ gives $P(\bar x+\bar a)+P(\bar x)=\bar b$, so $\bar x$ lies in a set of $\delta_P(\bar a,\bar b)\le2$
residues, since $P$ is APN. Each residue class modulo $2$ has $|2\GR(2^k,n)|=2^{n(k-1)}$ elements, so
$\delta^R_{\widehat P}(a,b)\le\delta_P(\bar a,\bar b)\,2^{n(k-1)}\le2\cdot2^{n(k-1)}$
(this is~\cite[Proposition~4(1)]{RS}).
\emph{Non-unit directions.} For $0\ne a\in(2)$, \cite[Proposition~5]{RS} gives
$\delta^R_{\widehat P}(a,b)\le t_{\max}\,2^{n(k-1)}$ with $t_{\max}$ the largest fibre of $\red(P)'$ on $F$,
which is $2$ because every nonempty fibre of $\red(P)'$ has exactly two elements. This is the only place where the
fibre condition is used.
\emph{Lower bound.} \cite[Theorem~3]{RS}: every permutation polynomial of $\GR(2^k,n)$ has
$\Delta^R\ge2\cdot2^{n(k-1)}$.
Hence $\Delta^R_{\widehat P}=2\cdot2^{n(k-1)}$. The last sentence of the corollary follows from the Main
Theorem.
\end{proof}

\begin{remark}[scope]\label{rem:ringscope}
Corollary~\ref{cor:ring} concerns permutation \emph{polynomials} of the ring; table-defined permutations of
$\GR(2^k,n)$ are not considered and can be better (\cite{RS} classifies APN permutations of $\GR(4,2)$).
Nothing is claimed about the security of any cipher. The theorem numbers of \cite{RS} used here refer to
its preprint version (IACR Cryptology ePrint Archive, 2026/1393).
\end{remark}

R{\o}njom and Sandrib already exhibit permutation polynomials over Galois rings with the optimal differential
uniformity $2\cdot2^{n(k-1)}$: the Dickson polynomials $g_5(x,u)=x^5-5ux^3+5u^2x$ \cite[Theorem~6]{RS}. Their residue
maps are not APN, so Corollary~\ref{cor:ring} is not implied by that example. For completeness we record the
short argument.

\begin{proposition}\label{prop:dickson}
Let $m\ge3$ be odd and $u\in\F_{2^m}^\times$. The map $D(x)=x^5+ux^3+u^2x$ of $\F_{2^m}$ (the reduction of
$g_5(x,u)$ modulo $2$) is not APN.
\end{proposition}
\begin{proof}
$D$ is quadratic, and for $a\ne0$ the map $\beta_a(x)=D(x+a)+D(x)+D(a)+D(0)=ax^4+a^4x+u(ax^2+a^2x)$ is
$\F_2$-linear; $D$ is APN iff $\ker\beta_a=\{0,a\}$ for every $a\ne0$. Put $x=ay$ and $z=y^2+y$. Since
$y^4+y=z^2+z$,
\[\beta_a(ay)=a^5(y^4+y)+ua^3(y^2+y)=a^3z\,\bigl(a^2z+a^2+u\bigr).\]
So $\beta_a(ay)=0$ iff $z=0$ or $z=c:=1+u/a^2$. The equation $y^2+y=c$ has two solutions iff $\Tr(c)=0$, and they
are different from $0$ and $1$ iff moreover $c\ne0$. Since $m$ is odd, $\Tr(1)=1$, so the condition is
$\Tr(u/a^2)=1$ with $u/a^2\ne1$. The elements of trace $1$ number $2^{m-1}\ge4$; choose one, $w\ne1$, and put
$a=\sqrt{u/w}$. Then $|\ker\beta_a|=4$ and $D$ is not APN.
\end{proof}

\noindent A direct computation over $\F_{2^m}$, $m\in\{3,5,7,9\}$, for every $u\ne0$ agrees with
Proposition~\ref{prop:dickson} (Supplement~S3).

\section{Formalization}\label{sec:lean}

The Main Theorem for odd $n\ge9$ is formalized in Lean~4~\cite{Lean4} with the Mathlib library~\cite{Mathlib}. The
final theorem is \lean{APNLift.LabelFormula.uopt\_odd}. In mathematical notation it states: for a finite field $F$
of characteristic $2$ with $|F|=2^n$, $n\ge9$ odd, and $Q:F\to F$ quadratic, APN and injective, there are additive
injective maps $L_1,A_1:F\to F$ such that $x\mapsto L_1(Q(A_1(x)))$ is an optimal representative. The definitions
used in the statement are the following.
\begin{itemize}
\item \lean{IsQuadratic Q}: for all $a,x,y$, $B(a,x+y)=B(a,x)+B(a,y)$, where $B(a,x)=Q(a+x)+Q(a)+Q(x)+Q(0)$
(Definition~\ref{def:quad}).
\item \lean{IsAPN Q}: for all $a\ne0$ and all $b$, $\#\{x:Q(x+a)+Q(x)=b\}\le2$ (Definition~\ref{def:apn}).
\item The reduced representative is the Lagrange interpolant of the map on all of $F$.
\item \lean{IsOptimalRepR38}\footnote{The suffix of this name is a version tag inside the formal development; the name is quoted as it occurs in the source.} \lean{h}: $h$ is quadratic, APN and bijective; the derivative of its reduced representative has
no zero on $F$; and there is $r\ne0$ with $\red(h)'(x)=\red(h)'(y)$ iff $y\in\{x,x+r\}$
(Definition~\ref{def:opt} and Remark~\ref{rem:optaffine}).
\end{itemize}
Additive injective maps of a finite field are exactly the invertible $\F_2$-linear maps, so the formal statement is
the Main Theorem for odd $n\ge9$. No hypothesis beyond those displayed remains: the existence of a good input, the
nonvanishing of $\Phi$ and $R$, and the linearization are proved inside the development.

The formal proof follows Sections~\ref{sec:layer1}--\ref{sec:mainproof}: the output constructions of
Section~\ref{ss:outputs}, the linearization of Theorem~\ref{thm:lin}, the Moore-point arguments of
Section~\ref{sec:cofactor}, and the count of Theorem~\ref{thm:count} with the degree $E_n=n(n-1)$ and the inequality of
Theorem~\ref{thm:universal} for all odd $n\ge9$. The Schwartz--Zippel bound is taken from Mathlib. A table
matching each statement of this paper with its formal declaration is in Supplement~S1.

The development was checked as follows. The source files were fixed by a SHA-256 digest of their contents. A build
from a new, empty directory completed without errors. The final theorem and its main lemmas depend only on the axioms
\lean{propext}, \lean{Classical.choice} and \lean{Quot.sound}, and no \lean{sorry} occurs. The statement of the
final theorem was restated in an audit file and checked against the declaration. This detects a declaration that
no longer proves the displayed type. It does not detect a change of the definitions that the type names, nor a
definitionally equal restatement; the meanings of those definitions were compared by the author with the list
above. The independent checker \lean{lean4checker}~\cite{lean4checker}, run in its \lean{--fresh} mode on
the module containing the final theorem, replayed every declaration that the theorem depends on, including those of
Mathlib and the Lean core library, and accepted them. The versions, the digests and the commands are listed in
Supplement~S1, which also gives a script that repeats these checks and fails if any of them fails.

\emph{Not formalized.} The Main Theorem for $n=5,7$ (Section~\ref{sec:small}); Appendix~\ref{app:pfaff} and the
numerical margins of Remark~\ref{rem:numbers}; Lemma~\ref{lem:grcrit} and the ring estimates used for
Corollary~\ref{cor:ring}; Proposition~\ref{prop:dickson}; Corollary~\ref{cor:proportion} and
Remark~\ref{rem:algorithm}; Proposition~\ref{prop:nthree} and the recount in Remark~\ref{rem:small}. Appendix~\ref{app:gold} is formalized separately in the
same development, except the hyperplane description of the image and the trace form of the zero criterion in Theorem~\ref{thm:sstrace} and the two points
named after Theorem~\ref{thm:gold} (Supplement~S1).

\emph{Data and code availability.} The Lean development (math-source digest
\texttt{409cc345b44d15399a3d33dbf81393ea79c61de777e0c2a9a214740c291ea37a}; Supplement~S1), the certificates and input
witnesses for $n=5,7$, the scripts for Proposition~\ref{prop:nthree} and Corollary~\ref{cor:proportion}, and the
sources of this paper, its supplement and its blueprint are available at \url{https://github.com/hirokifukui/apnlift}
(release v1.0.3) and are archived at Zenodo under \url{https://doi.org/10.5281/zenodo.23053111} (all versions of the
source archive); the paper and its supplementary files are archived under
\url{https://doi.org/10.5281/zenodo.23053113}.

\section{Concluding remarks}\label{sec:concl}

Every quadratic APN permutation in odd dimension $n\ge5$ has an optimal representative $L\circ Q\circ A$ with $L,A$
invertible $\F_2$-linear. For $n\ge9$ the input $A$ comes from a pigeonhole comparison between the invertible binary
matrices on one side and, on the other, the matrices of rank at most $n-2$ together with a hypersurface of degree
$n(n-1)$; this case is also checked in Lean. For $n=5,7$ the input comes from the classifications of quadratic APN
functions and explicit witnesses, rechecked on a second classification. The range is exact: the statement fails
for $n=3$ and is empty for even $n$. For $n\ge9$ at least a fixed fraction of all invertible inputs works (more
than half for $n\ge15$), so an optimal representative can be found by independent uniform sampling of inputs.

Several questions remain open and are not needed for the results above.
\begin{itemize}
\item \emph{Higher degree.} For even $n$ there are no quadratic APN permutations~\cite{BCCL}, and $n=3$ is a genuine
exception (Proposition~\ref{prop:nthree}; see also~\cite{Stone-cx}). The argument uses quadratic $Q$ through the linearization of
Section~\ref{sec:layer2}; nothing is claimed for APN permutations of higher algebraic degree. Experiments on
relabelled power maps~\cite{Stone-cx}, including Welch functions of degree three, suggest that reduced representatives with
zero-free derivatives are common there as well.
\item \emph{Zero-derivative normal form.} Is every quadratic APN permutation in odd dimension $n\ge5$ linearly equivalent to
a permutation $S$ with $\red(S)'=0$? For such $S$, Proposition~\ref{prop:outNF} gives an optimal map $S\circ A_t$
directly, without the counting argument. This holds for the Gold classes (Appendix~\ref{app:gold}) and for further families
(Supplement~S4).
\item \emph{A uniform argument for $n=5,7$.} The counting bound of Theorem~\ref{thm:count} fails for $n=5,7$; an
argument that avoids the classifications would also make these cases amenable to formalization.
\item \emph{Same input.} Theorem~\ref{thm:twoinv} is a sufficient condition at a given input. An optimal output can
also exist at an input that fails the criterion, but not always: for $n=9$ and $n=11$ there are Dembowski--Ostrom (purely quadratic: every exponent of $\red(P)$ has binary weight
exactly $2$) APN permutations $P$
with $s=1$ and $p_P=0$ at the identity input; then $\red(L\circ P)'=K_\ell$ vanishes at $0$ for every output $L=L_\ell$,
so no output at that input is optimal (Supplement~S2). A description of the inputs at which an optimal output exists
is not known.
\end{itemize}

\appendices
\section{A sharper degree: the Pfaffian route}\label{app:pfaff}

One can replace the adjugate entry $R$ of Proposition~\ref{prop:R} by a polynomial $c$ built from Pfaffians, of
degree $\tfrac{(n-1)(n-2)}2-1$. This lowers the degree of the exceptional hypersurface from $E_n=n(n-1)$ to
$D_n=\tfrac{n^2+n-4}{2}$ and increases the margins (Remark~\ref{rem:numbers}), but it does not change what is
proved: both routes give Theorem~\ref{thm:universal} for all odd $n\ge9$ and neither closes at $n=5,7$. This
variant is not formalized.

\begin{lemma}[Pfaffian vector]\label{lem:pf}
Let $M$ be an alternating $n\times n$ matrix over a field of characteristic $2$, $n$ odd, and
$\pi(M)_i=\Pf(M_{\hat i\hat i})$, the Pfaffian of the principal submatrix without row and column $i$. Then
$M\pi(M)=0$, and $\pi(M)\ne0$ iff $\rank M=n-1$. Each coordinate of $\pi$ is homogeneous of degree $(n-1)/2$
in the entries of $M$.
\end{lemma}
\begin{proof}
Fix $r$. Border $M$ by a row and a column with index $0$ equal to row $r$ of $M$:
$\tilde M_{0,j+1}=\tilde M_{j+1,0}=M_{r,j}$, $\tilde M_{i+1,j+1}=M_{ij}$. Then $\tilde M$ is alternating of
even order $n+1$, and rows $0$ and $r+1$ coincide (both are $(0,M_{r,\cdot})$, since $M_{rr}=0$), so
$\Pf(\tilde M)^2=\det\tilde M=0$. Expanding $\Pf(\tilde M)$ along row $0$ gives
$\sum_jM_{rj}\Pf(M_{\hat j\hat j})=(M\pi(M))_r$ (no signs in characteristic $2$).
Alternating matrices have even rank, so $\rank M\le n-1$. If $\rank M=n-1$, let $r\ne0$ span the radical and
choose $i$ with $r_i\ne0$; the coordinate hyperplane $\{x_i=0\}$ is a complement of the radical, so the form is
nondegenerate on it, $\det M_{\hat i\hat i}\ne0$ and $\Pf(M_{\hat i\hat i})\ne0$. If $\rank M\le n-3$, every
principal $(n-1)$-submatrix is singular, so $\pi(M)=0$.
\end{proof}

Put $P(Z):=\pi(\bG_{\omega(Z)})\in F[Z]^n$; each coordinate is zero or homogeneous of degree
$\tfrac{n-1}{2}(n-2)$.

\begin{lemma}[the polynomial $c$]\label{lem:cpf}
Let $Q$ be a quadratic APN permutation, $n\ge5$ odd. There is a nonzero homogeneous $c\in F[Z]$ of degree
$\frac{(n-1)(n-2)}2-1$ with $P(Z)=c(Z)\,Z$. For $\zeta\in F^n\setminus0$:
$c(\zeta)\ne0$ iff $\omega(\zeta)\ne0$ and $\rank\bG_{\omega(\zeta)}=n-1$; in that case
$\ker\bG_{\omega(\zeta)}=F\zeta$.
\end{lemma}
\begin{proof}
Work over the algebraic closure $\bar F$; Theorem~\ref{thm:lin}(c), \eqref{eq:GGbar} and Lemmas~\ref{lem:s}, \ref{lem:pf} are polynomial
identities or linear algebra valid over any extension field. Let $O=\{\zeta\in\bar F^n:P(\zeta)\ne0\}$. For
$\zeta\in O$ we have $\bG_{\omega(\zeta)}\ne0$, so $\omega(\zeta)\ne0$ (and $\zeta\ne0$, since $P(0)=0$). By
Lemma~\ref{lem:s}, $G(\zeta)\omega(\zeta)=0$, i.e.\ $\bG_{\omega(\zeta)}\zeta=0$ by \eqref{eq:GGbar}. By
Lemma~\ref{lem:pf}, $\rank\bG_{\omega(\zeta)}=n-1$, and both $\zeta$ and $P(\zeta)\ne0$ lie in its
one-dimensional kernel. Hence $P_i(\zeta)\zeta_j=P_j(\zeta)\zeta_i$ on $O$.

\emph{$O$ is not empty.} Take $x\in F^\times$ and $\zeta=m(x)$. By \eqref{eq:GGbar} and (MB),
$\ker G(m(x))=\{\ell:\Gamma(m(x),m(z))\cdot\ell=0\ \forall z\in F\}=\{\ell:L_\ell\text{ vanishes on }
H_x:=B_Q(x,F)\}$. APN gives $\ker B_Q(x,\cdot)=\{0,x\}$, so $H_x$ is an $\F_2$-hyperplane
$\{z:\Tr(a_xz)=0\}$ with $a_x\ne0$. The $\F_2$-linear maps vanishing on $H_x$ are $y\mapsto b\,\Tr(a_xy)$,
$b\in F$, i.e.\ $L_\ell$ with $\ell\in F\,m(a_x)$. So $\rank G(m(x))=n-1$ and, by Lemma~\ref{lem:s},
$\omega(m(x))=t\,m(a_x)$ with $t\ne0$. For $a\ne0$, $\bG_{m(a)}$ is the Gram matrix of the $F$-bilinear form
$(\xi,\chi)\mapsto\Gamma(\xi,\chi)\cdot m(a)$, which on Moore vectors equals
$L_{m(a)}(B_Q(x',z))=\Tr(aB_Q(x',z))=b_a(x',z)$. In the $F$-basis $m(x_1),\dots,m(x_n)$ of (MB) its Gram
matrix is the $\F_2$-Gram matrix of $b_a$, whose rank is $n-1$ by (AB); rank does not change under congruence
or field extension. So $\rank\bG_{\omega(m(x))}=n-1$, and $P(m(x))\ne0$ by Lemma~\ref{lem:pf}.

$O$ is a nonempty Zariski-open subset of the irreducible $\bar F^n$, hence dense, so
$P_iZ_j-P_jZ_i\equiv0$. Then $Z_0\mid P_0Z_1$ gives $Z_0\mid P_0$; with $c:=P_0/Z_0$, $P_iZ_0=P_0Z_i=cZ_0Z_i$
gives $P_i=cZ_i$. Since $P\ne0$, $c\ne0$, and $c$ is homogeneous of degree $\deg P-1$. For
$\zeta\in F^n\setminus0$: $c(\zeta)\ne0$ iff $P(\zeta)\ne0$, which by Lemma~\ref{lem:pf} holds iff
$\rank\bG_{\omega(\zeta)}=n-1$ (this forces $\omega(\zeta)\ne0$). The kernel statement was shown above.
\end{proof}

\begin{lemma}[injectivity, Pfaffian form]\label{lem:injpf}
The map $[\zeta]\mapsto[\omega(\zeta)]$ is injective on $\{[\zeta]\in\mathbb P^{n-1}(F):c(\zeta)\ne0\}$.
\end{lemma}
\begin{proof}
As Lemma~\ref{lem:inj}, with $\ker\bG_{\omega(\zeta)}=F\zeta$ from Lemma~\ref{lem:cpf}.
\end{proof}

\begin{proposition}[counting with $D_n$]\label{prop:countD}
For every quadratic APN permutation $Q$ of $F$, $n$ odd, $n\ge5$,
$|\mathrm{Fail}(Q)|\le\rho_n(\le n-2)-1+D_nq^{n-1}$, and for odd $n\ge9$,
$|\GL_n(\F_2)|-|\mathrm{Fail}(Q)|\ge(\Lambda_n-D_n2^{-n})q^n+1$.
\end{proposition}
\begin{proof}
The proof of Theorem~\ref{thm:count} with $R$ replaced by $c$ (Lemma~\ref{lem:cpf}, Lemma~\ref{lem:injpf}) and
$\deg\Phi c\le(2n-2)+\tfrac{(n-1)(n-2)}2-1=D_n$. Positivity follows from Theorem~\ref{thm:universal}, since
$D_n\le E_n$.
\end{proof}

\section{An explicit family for the Gold classes}\label{app:gold}

This appendix is not used in the proof of the Main Theorem. It gives closed-form optimal representatives in every
Gold class, with no estimate and no search, by an input shear.

For $R:F\to F$ write $H_R(x)=R(x+1)+R(x)+R(1)+R(0)$ and $s_R=R(1)+R(0)$; for quadratic $R$ the map
$H_R$ is additive. For $t\in F$ put $A_t(x)=x+\Tr(tx)$; if $\Tr(t)=0$ then $A_t$ is an additive involution,
and $T_t=\sum_{i<n}t^{2^i}X^{2^i}$ represents $x\mapsto\Tr(tx)$. For quadratic $f$ every exponent of
$\red(f)$ has binary weight at most $2$.

\subsection{The square-shear criterion}
\begin{theorem}[square shear]\label{thm:ss}
Let $n\ge2$ and let $R$ be quadratic with $\red(H_R)$ supported on $\{X^{2^k}:1\le k\le n-1\}$; let $h$ be
its coefficient of $X^{2^{n-1}}$. Then
$\red(R\circ A_t)=\red(R)+\bigl(\red(H_R)+s_R\bigr)T_t-h\,t^{q/2}\,(X^q-X)$. If moreover $\red(R)'=0$, then
$\red(R\circ A_t)'=(\red(H_R)+s_R)\,t+h\,t^{q/2}$.
\end{theorem}
\begin{proof}
$R(A_t(x))=R(x)+\Tr(tx)(H_R(x)+s_R)$: if $\Tr(tx)=0$ both sides are $R(x)$, otherwise
$R(x+1)=R(x)+H_R(x)+s_R$. So the right-hand side represents $R\circ A_t$, since $X^q-X$ vanishes on $F$.
The product $(\red(H_R)+s_R)T_t$ has degree at most $q$ and coefficient $h\,t^{2^{n-1}}$ at $X^q$, which
the correction removes; so the right-hand side has degree $<q$. $\red(H_R)$ has only even exponents,
$T_t'=t$ and $(X^q-X)'=-1$.
\end{proof}

\begin{theorem}[trace form]\label{thm:sstrace}
Let $R$ be a quadratic APN permutation with $\red(R)'=0$ and $\red(H_R)$ supported on $X^{2^k}$,
$1\le k\le n-1$, with coefficient $h$ at $X^{2^{n-1}}$. There is $\pi\ne0$ with
$\{R(x+1)+R(x):x\in F\}=\{y:\Tr(\pi y)=1\}$. For $t\ne0$ with $\Tr(t)=0$, $R\circ A_t$ is an APN
permutation, $\red(R\circ A_t)'(x)=t\,(R(x+1)+R(x))+h\,t^{q/2}$ is exactly two-to-one with fibres
$\{x,x+1\}$, and it has no zero on $F$ iff $\Tr(\pi^2h^2/t)=0$.
\end{theorem}
\begin{proof}
$H_R$ is $\F_2$-linear with kernel $\{0,1\}$ (APN), so its image $V$ is a hyperplane $\{\Tr(\pi y)=0\}$;
since $R$ is injective, $s_R\notin V$ and $\{H_R(x)+s_R\}=V+s_R$. $R\circ A_t$ is a quadratic bijection whose
polar kernels are images of those of $R$, so it is APN. By Theorem~\ref{thm:ss} the derivative is
$t(R(x+1)+R(x))+ht^{q/2}$; two values coincide iff $x+y\in\ker H_R=\{0,1\}$. A zero exists iff
$ht^{q/2}/t$ lies in $V+s_R$, i.e.\ iff $\Tr(\pi ht^{q/2-1})=1$; and $(\pi ht^{q/2-1})^2=\pi^2h^2/t$.
\end{proof}

\begin{theorem}\label{thm:gold}
Let $n\ge5$ be odd, $1\le r<n$ with $\gcd(r,n)=1$, $e=\min(r,n-r)$, and $t\in F^\times$ with $\Tr(t)=0$.
Put $T_t(X)=\sum_{i<n}t^{2^i}X^{2^i}$ and $P_{e,t}=X^{2^{e+1}+2}+(X^{2^{e+1}}+X^2+1)\,T_t$. Then $P_{e,t}$
is an optimal APN permutation, linearly equivalent to $x^{2^r+1}$, and $P_{e,t}'=t(X^{2^{e+1}}+X^2+1)$. The
$q/2-1$ maps $P_{e,t}$ are distinct.
\end{theorem}

\begin{proof}
Since $x^{2^{n-r}+1}=(x^{2^r+1})^{2^{n-r}}$, the classes of $r$ and $n-r$ coincide, so take
$r=e\le(n-1)/2$. Let $R(x)=x^{2^{e+1}+2}$, an output Frobenius twist of the Gold function; it is a quadratic
APN permutation with $\red(R)'=0$. Then $H_R(x)=x^{2^{e+1}}+x^2$, $s_R=1$, and $h=0$ since $2\le e+1\le n-2$.
Theorem~\ref{thm:ss} gives $\red(R\circ A_t)=P_{e,t}$ and $P_{e,t}'=t(X^{2^{e+1}}+X^2+1)$. The value
$x^{2^{e+1}}+x^2+1$ vanishes iff $R(x+1)=R(x)$, which is impossible. Fibres are as in
Theorem~\ref{thm:sstrace}. $P_{e,t}$ determines $T_t$ and hence $t$, so the $q/2-1$ maps are distinct.
\end{proof}
For $r=e$ the assertions of Theorem~\ref{thm:gold} (bijectivity, the APN property, the reduced polynomial and its
derivative, absence of zeros, the fibres, distinctness in $t$) are formalized in the development of
Section~\ref{sec:lean}; so is the identity that reduces $r>n/2$ to $r=e$. The linear equivalence is the construction
$x\mapsto\bigl(A_t(x)^{2^e+1}\bigr)^2$ itself and is not a separate formal statement (Supplement~S1). The theorem was
also checked numerically for $1009$ parameter sets (Supplement~S4). The explicit counterexamples
of~\cite{Stone-cx} are of the same kind, a Gold power map composed with an $\F_2$-linear input change (for instance
$(x+x^8+x^{16})^3$ on $\F_{2^5}$); Theorem~\ref{thm:gold} gives such maps in closed form for every odd $n\ge5$ and
every Gold exponent.

\section*{Acknowledgment}
\emph{Use of AI systems.} In accordance with the IEEE policy on AI-generated content, the author discloses the
following substantive uses. Claude (Anthropic), including the agentic tool Claude Code, and ChatGPT (OpenAI) were
used as follows.
\begin{itemize}
\item \emph{Mathematics (Sections~\ref{sec:layer1}--\ref{sec:cons}, Appendices~\ref{app:pfaff}--\ref{app:gold}).}
Proof exploration and drafting of arguments: AI systems proposed and drafted proofs, which were then revised; ChatGPT
was also used for repeated adversarial reviews of the whole argument.
\item \emph{Formalization (Section~\ref{sec:lean}, Supplement~S1).} The Lean code was written with Claude and
Claude Code; its correctness rests on the Lean kernel and the independent checker, not on the AI systems.
\item \emph{Computations (Section~\ref{sec:small}, Supplements~S2--S4).} Scripts for the finite computations,
certificates and witnesses were designed, written and run with AI assistance; the results were repeated with
independently written code and, where stated, with a second computer algebra system.
\item \emph{Text.} The exposition of all sections was drafted and edited with AI assistance.
\end{itemize}
No output of these systems was accepted without checking: the statements listed in Section~\ref{sec:lean} were checked
by the Lean kernel; the computations were repeated as stated; the remaining arguments were checked by the author and
by further AI-assisted reading passes. These passes are internal checks, not external peer review. The author takes
full responsibility for the definitions, statements, proofs, citations and accompanying material.

\bibliographystyle{IEEEtran}
\bibliography{references}
\end{document}